\documentclass{amsart}
\usepackage{amssymb}
\usepackage{stmaryrd}
\usepackage[all]{xy}
\usepackage{epsf}

\newtheorem{theorem}{Theorem}
\newtheorem{proposition}{Proposition}
\newtheorem{lemma}{Lemma}
\newtheorem{corollary}{Corollary}

\newtheorem{descr}{Description}
\newtheorem{definition}{Definition}

\theoremstyle{definition}
\newtheorem{example}{Example}

\theoremstyle{remark}
\newtheorem{remark}{Remark}

\newcommand{\Rr}{\mathbb R}
\newcommand{\Zz}{\mathbb Z}
\newcommand{\Cc}{\mathbb C}

\newcommand{\set}[1]{\left\{#1\right\}}

\newcommand{\eps}{\varepsilon}

\newcommand{\rmap}{\longrightarrow}

\newcommand{\X}{\ensuremath{\mathfrak{X}}}
\newcommand{\F}{\ensuremath{\mathcal{F}}}
\renewcommand{\P}{\ensuremath{\mathcal{P}}}
\renewcommand{\S}{\ensuremath{\mathcal{S}}}

\newcommand{\NN}{\ensuremath{\mathcal{N}}}
\newcommand{\AAA}{\ensuremath{\mathcal{A}}}

\newcommand{\G}{\mathcal{G}}            % Lie groupoid
\newcommand{\s}{\mathbf{s}}             % source map
\renewcommand{\t}{\mathbf{t}}           % target map
\renewcommand{\H}{\mathcal{H}}          % Lie subgroupoid
\DeclareMathOperator{\Ad}{ad}           % Adjoint
\DeclareMathOperator{\Ker}{Ker}         % Kernel
\DeclareMathOperator{\rank}{rank}       % rank of a vector bundle
\DeclareMathOperator{\codim}{codim}     % codimension
\DeclareMathOperator{\Ima}{Im}          %
\renewcommand{\Im}{\Ima}                % Image

\newcommand{\al}{\alpha}                % section of Lie algebroid
\newcommand{\be}{\beta}                 % section of Lie algebroid
\newcommand{\Lie}{\mathcal{L}}          % Lie derivative
\renewcommand{\gg}{\mathfrak{g}}        % Lie algebra

\newcommand{\var}{var}                  % variation

\begin{document}

\title{Integrability of Poisson brackets}
%
% author one information
\author{Marius Crainic}
\address{Depart. of Math., Utrecht University, 3508 TA Utrecht, 
The Netherlands}
\email{crainic@math.uu.nl}
\thanks{Supported in part by NWO and a Miller Research Fellowship}

% author two information
\author{Rui Loja Fernandes}
\address{Depart.~de Matem\'{a}tica, 
Instituto Superior T\'{e}cnico, 1049-001 Lisboa, PORTUGAL} 
\email{rfern@math.ist.utl.pt}
\thanks{Supported in part by FCT/POCTI/FEDER and by grant
  POCTI/1999/MAT/33081.}

%%% ----------------------------------------------------------------------
\begin{abstract}
We discuss the integration of Poisson brackets, motivated by our
recent solution to the integrability problem for general Lie
brackets. We give the precise obstructions to integrating Poisson
manifolds, describing the integration as a symplectic quotient, in the
spirit of the Poisson sigma-model of Cattaneo and Felder. For regular
Poisson manifolds we express the obstructions in terms of variations
of symplectic areas, improving on results of Alcade Cuesta and Hector.
We apply our results (and our point of view) to decide about the
existence of complete symplectic realizations, to the integrability of
submanifolds of Poisson manifolds, and to the study of dual pairs,
Morita equivalence and reduction.
\end{abstract}
%%% ----------------------------------------------------------------------

\maketitle

\tableofcontents

\newpage
%%%%%%%%%%%%%%%%%%%%%%%%%%%%%%%%%%%%
%%%%%%%%%%%%%%%%%%%%%%%%%%%%%%%%%%%%
%%%%%%%%%%%%%%%%%%%%%%%%%%%%%%%%%%%%
\section{Introduction}             %
\label{Integrability: Outline}     %
%%%%%%%%%%%%%%%%%%%%%%%%%%%%%%%%%%%%
%%%%%%%%%%%%%%%%%%%%%%%%%%%%%%%%%%%%
%%%%%%%%%%%%%%%%%%%%%%%%%%%%%%%%%%%%

A \textbf{Poisson bracket} on a manifold $M$ is a Lie bracket $\{\cdot,
\cdot\}$ on the space $C^{\infty}(M)$ of smooth functions on $M$,
satisfying the derivation property
\[ \{ fg, h\}= f\{ g, h\}+ g\{ f, h\}, \quad f,g,h \in
C^{\infty}(M).\] The integrability problem for Poisson brackets can be
loosely stated as: 
\begin{itemize}
\item Is there a Lie group integrating this Lie algebra?
\end{itemize}
Stated as such, this problem is beyond our current state of knowledge,
as is illustrated by the well known \emph{flux conjecture} in
symplectic geometry \cite{MaSa} (see also Milnor's remarks about
infinite dimensional groups in \cite{Mil}). However, such problems
become more tractable if instead of infinite dimensional Lie groups
one brings finite dimensional objects known as \emph{Lie groupoids}
into the picture. 

A Poisson bracket on a manifold $M$ gives rise to a Lie bracket
$[\cdot, \cdot]$ on the space $\Omega^1(M)$ of 1-forms on $M$. This
bracket is uniquely determined by the following two requirements:
\begin{enumerate}
\item[(i)] for exact 1-forms it coincides with the Poisson bracket: 
\[ [df, dg]= d\{f, g\},\quad f,g\in C^{\infty}(M);\] 
\item[(ii)] it satisfies the Leibniz identity:
\begin{equation}
\label{Leibniz} [\al, f\be]= f[\al, \be] + \#\be(f)\al, \quad
\al,\be\in\Omega^1(M),\ f\in C^{\infty}(M).
\end{equation}
\end{enumerate}
Here $\#:T^*M\to TM$ denotes contraction by the Poisson 2-tensor
$\Pi\in \Gamma(\Lambda^2TM)$ which is associated to the Poisson
bracket by $\Pi(df, dg)= \{f, g\}$. An explicit
formula for this bracket is
\begin{equation}
\label{kozul:bracket} [\al, \be]=
\Lie_{\#\al}\be-\Lie_{\#\be}\al-d\Pi(\al,\be),\quad
\al,\be\in\Omega^1(M).
\end{equation}
The triple $(T^*M,[\cdot, \cdot],\#)$ is an example of a \emph{Lie
algebroid}. Lie algebroids can be thought as 
``infinite dimensional Lie algebras of geometric type'', or 
``generalized tangent bundles''.
In general, a \textbf{Lie algebroid} is a vector bundle $A\to
M$ together with a Lie bracket $[\cdot, \cdot]$ on the space of
sections $\Gamma(A)$ and a bundle map $\#:A\to TM$, giving rise to a
Lie algebra morphism $\#:\Gamma(A)\to\X(M)$, such that the
analogue of Leibniz's identity (\ref{Leibniz}) is satisfied for all
$\al,\be \in \Gamma(A)$.

The global counterpart to Lie algebroids are \emph{Lie groupoids}. A
Lie groupoid consists of arrows (transformations) between different
objects (points), which can be (smoothly) multiplied provided they
match. More formally, a groupoid is a small category where every
morphism is an isomorphism. A \textbf{Lie groupoid} is a groupoid in
the differentiable category: it consists of a manifold $\Sigma$ (the
arrows), together with two submersions $\s, \t: \Sigma \to M$ (the
source and the target maps) onto the base manifold $M$ (the objects),
an embedding $M\to \Sigma$, $x\mapsto 1_{x}$ (the unit section), and a
smooth map $\Sigma\times_{M}\Sigma\to \Sigma$, $(g, h)\mapsto gh$ (the
multiplication) defined on the space of pairs $(g, h)$ with $\s(g)=
\t(h)$. We will follow the conventions of \cite{CaWe}.

To every Lie groupoid there is associated a Lie algebroid (see
\cite{CaWe}). The converse is not true, and the precise obstructions
to the integration of Lie algebroids to Lie groupoids were determined in
\cite{CrFe}.  Uniqueness, up to isomorphism, can be obtained by
requiring the $\s$-fibers to be simply connected \footnote{in this paper, 
by simply connected we mean connected with vanishing fundamental group}
and, given $\Sigma$, one can construct a $\s$-simply connected groupoid
$\widetilde{\Sigma}$ by taking universal coverings of the $\s$-fibers
(see \cite{MoMr,CrFe}).  The integrability problem for Poisson manifolds 
can then be restated as:
\begin{itemize}
\item Is there a Lie groupoid integrating $T^*M$?
\end{itemize}
Note that the integrability of $T^*M$ really
amounts to the integrability of the brackets $[\cdot, \cdot]$ on
$\Omega^1(M)$.

The integration of Poisson manifolds is, of course, a particular case
of the integrability problem for general Lie algebroids. As already
mentioned above, a complete
solution for the latter was presented in \cite{CrFe}.  However, it is
important to understand this special case. On the one hand, the
integration of Poisson manifolds is richer due to the presence of
symplectic geometry on the leaves
of the characteristic foliation,
which gives rise to many properties that are not present in the
general integrability problem. On the other hand, the integrability
problem for Poisson brackets is relevant, for example, to symplectic
reduction \cite{Daz2}, to Poisson topology \cite{Ginz}, to various
quantization schemes \cite{CaFe,Ko,Wein}, and to many other
problems. As an example of the richer geometry presented in the
special case of Poisson manifolds, we shall see that a Lie groupoid
integrating a Poisson manifold has a natural symplectic structure, a
well-known fact going back to the earlier works of Weinstein \emph{et al}.

Recall that a \textbf{symplectic groupoid} is a Lie groupoid $\Sigma$
together with a symplectic form, such that the graph of the
multiplication is Lagrangian. This apparently mild condition is
actually quite strong. It induces a natural Poisson structure on $M$,
satisfying the following properties (see \cite{CDW}):
\begin{enumerate}
\item[(a)] $\s$ is Poisson and $\t$ is anti-Poisson(\footnote{A
\emph{Poisson map} is a map $\phi:M\to N$ between two Poisson
manifolds that preserves the Poisson brackets. A Poisson map is called
\emph{complete} if, whenever $X_h$ is a complete Hamiltonian vector
field on $N$, then $X_{\phi^*h}$ is also a complete Hamiltonian vector
fields on $M$.});
\item[(b)] both $\s$ and $\t$ are complete maps;
\item[(c)] the $\s$-fibers and the $\t$-fibers are symplectic orthogonal;
\item[(d)] $M$, viewed as the unit section, is a Lagrangian submanifold
of $\Sigma$;
\item[(e)] the Lie algebroid of $\Sigma$ is canonically
isomorphic to $T^*M$.
\end{enumerate}
The inverse problem (reconstructing $\Sigma$) is the symplectic
version of the integrability problem:
\begin{itemize}
\item Is there a symplectic groupoid integrating $M$?
\end{itemize}
Again, it is not hard to see that if a Poisson manifold admits an
integrating symplectic groupoid, then it admits a unique $\s$-simply
connected one. In this case we say that $M$ is an integrable Poisson
manifold, and the $\s$-simply connected groupoid $\Sigma=\Sigma(M)$
integrating $M$ will be called the symplectic groupoid of $M$.

By property (e) above, if a Poisson manifold is integrable, then
its associated algebroid $T^*M$ is also integrable. The converse
is also true as was shown by Mackenzie and Xu in \cite{MaXu}:

\begin{theorem}
If $M$ is a Poisson manifold such that $T^*M$ is an integrable Lie
algebroid then $M$ is an integrable Poisson manifold.
\end{theorem}

Many known results like this have a different and simpler proof in our
unified approach to the integrability problem to be presented below.
Let us give a short overview of our solution to the integrability
problem and, at the same time, an outline of the content of the paper.

In \cite{CrFe}, for \emph{any} Lie algebroid $A$ we have constructed a
topological groupoid $\G(A)$, called the Weinstein groupoid of $A$
which is a fundamental invariant of the Lie algebroid. Moreover,
this groupoid has a compatible differentiable structure (i.e., is a
Lie groupoid) if and only if $A$ is an integrable Lie algebroid. In the case of
Poisson manifolds, where $A=T^*M$, we will denote $\G(T^*M)$ by
$\Sigma(M)$. One should think of $\Sigma(M)$ as the homotopy
(or fundamental) 
groupoid 
of the Poisson manifold, and in fact it can be described,
as we shall explain in Section \ref{Cotangent-paths} below, as a
quotient
\[ \Sigma(M)=\text{cotangent paths}/\text{cotangent homotopies}.\]
The obvious similarity to the ordinary homotopy group is related with
the following basic \emph{philosophical principle}:
\begin{quote}
\emph{In analogy with the role played by the tangent bundle $TM$ of
manifolds, Lie algebroids can be thought of as the ``tangent bundles of
singular structures''.}
\end{quote}
Accordingly, many constructions/results in differential geometry that
can be carried out in terms of just the tangent bundle make sense for
general Lie algebroids. However, sometimes this is not entirely
obvious, and the Lie algebroid framework can reveal totally new
patterns. For instance, in contrast with the fact that the fundamental groupoid of a
manifold is always smooth, $G(A)$ may fail to be smooth. The main result of \cite{CrFe} gives the precise
obstructions for a differentiable structure to exist in $\G(A)$, and
expresses them in terms of so-called \emph{monodromy groups}. For the
special case of Poisson manifolds, we shall describe
them in Section
\ref{monodromy} below.

There are two special properties of $T^*M$ that distinguish this case
from the general case. First, the anchor and the bracket are both
induced from the Poisson tensor and so are intimately related. Second,
one has the duality between $T^*M$ and $TM$. As we shall see in
Section \ref{symplectic-groupoid}, these lead to a description of
$\Sigma(M)$ as a symplectic quotient, a fact first noted in
\cite{CaFe}. Then the symplectic form on $\Sigma(M)$, makes the
Weinstein groupoid into a symplectic Lie groupoid, hence reproving
Mackenzie-Xu's result mentioned above.

Our main result, whose proof is given in Section \ref{smoothness}
below, can be stated as follows:

\begin{theorem}\label{main-theorem} For a Poisson manifold $M$, the
following are equivalent:
\begin{enumerate}
\item[(i)] $M$ is integrable by a symplectic Lie groupoid,
\item[(ii)] the algebroid $T^*M$ is integrable,
\item[(iii)] the Weinstein groupoid $\Sigma(M)$ is a smooth manifold,
\item[(iv)] the monodromy groups $\NN_{x}$, with $x\in M$, are locally uniformly
discrete;
\end{enumerate}
\end{theorem}

The last condition of the theorem makes no reference to
algebroids/groupoids, and gives an integrability criteria which is
computable in explicit examples. Moreover, the monodromy groups
$\NN_{x}$ (which are just some additive subgroups of the co-normal
vector space $\nu_{x}^{*}(L)$ to the symplectic leaf $L$ through $x$)
are invariants of Poisson manifolds which are interesting on their
own. For regular Poisson manifolds, as we discuss in Section
\ref{Integrability:regular} (improving on results of Alcade
Cuesta-Hector \cite{AlHe}), the monodromy can be expressed in terms of
variations of symplectic areas of spheres along transverse directions
to the symplectic leaves. In Section \ref{Integrability:general},
we present many examples of integrable and non-integrable Poisson
manifolds.

A well-known result of Karasev and Weinstein states that \emph{any}
Poisson manifold $M$ has a \emph{symplectic realization}, i.e, there
exists a symplectic manifold $S$ and a surjective Poisson submersion
$\mu:S\to M$. The existence of symplectic realizations with $\mu$
\emph{complete} has been an open problem, which can be thought as yet another
instance of the integrability problem for
  Poisson manifolds:
\begin{itemize}
\item Is there a complete symplectic realization of $M$?
\end{itemize}
In Section \ref{realizations} below we solve this problem 
(and sketch a different proof of the result of Karasev and Weinstein):

\begin{theorem}
A Poisson manifold admits a complete symplectic realization if and
only if it is integrable. 
\end{theorem}

This establishes the equivalence of all the different notions of
integrability. In the last two sections of the paper we give some
more applications of our integrability results. 

In Section \ref{Poisson:category} we consider the Poisson brackets
induced on a submanifold of a Poisson manifold. 
We clarify
and improve the results of Xu \cite{Xu1} and Vaisman \cite{Vai2} on
induced Poisson structures. We also discuss, for an integrable Poisson
manifold, when is the induced Poisson bracket on a submanifold
integrable. 

In Section \ref{Morita} we discuss Morita equivalence of Poisson
manifolds, both for integrable and non-integrable Poisson
manifolds. We observe that the original definition due to Xu only makes
sense for integrable Poisson manifolds. For non-integrable Poisson
manifolds we consider a weak notion of Morita equivalence, and we
prove that many invariants of Poisson manifolds, relevant for the
integrability problem, are \emph{weak} Morita invariant. 

\begin{remark}[Hausdorff Issues] 
\label{rem:Hausdorff}
The reader will notice that we are forced to allow non-Hausdorff
manifolds in our paper. There are at least two simple reasons for
this. First of all, a bundle of Lie algebras $\gg$ over a manifold
$B$
may not integrate to a bundle 
$G$ of Lie groups over $B$ (that
is, the Lie algebra of the fiber Lie group $G_{b}$ coincides with
$\gg_{b}$, for all $b\in B$), if we require $G$ to be
Hausdorff. However, there exists always at least one which is a
(possibly non-Hausdorff) manifold (see \cite{DoLa}). Secondly, there
are simple examples of foliations whose graph, although a manifold,
may be non-Hausdorff.

We mention here the objects of this paper which allow non-Hausdorff manifolds:
\begin{enumerate}
\item[(a)] Lie groupoids (Sections \ref{symplectic-groupoid} and
\ref{smoothness}) may be non-Hausdorff manifolds. However, the base
space (denoted $M$ here), as well as the $\s$-fibers and $\t$-fibers,
are always assumed to be Hausdorff.
\item[(b)] For a symplectic realization $\pi: S\to M$ (Section
\ref{realizations}), $S$ is allowed to be non-Hausdorff. However, the
leaves of the foliation $\F(\pi)$ by fibers of $\pi$ and of the
symplectic orthogonal foliation $\F(\pi)^{\perp}$, are all Hausdorff.
\end{enumerate}
In particular, the Poisson manifolds we study are always assumed to be
Hausdorff and paracompact. 
\end{remark}

\begin{center}
\textsc{Acknowledgments}
\end{center}

During the preparation of this manuscript we were influenced by
discussions with many colleagues. We would especially
like to thank
Henrique Bursztyn and Alan Weinstein for comments and suggestions.

%%%%%%%%%%%%%%%%%%%%%%%%%%%%%%%%%%%%%%%%%%%%%%%%%%%%%%%%%%%%
%%%%%%%%%%%%%%%%%%%%%%%%%%%%%%%%%%%%%%%%%%%%%%%%%%%%%%%%%%%%
%%%%%%%%%%%%%%%%%%%%%%%%%%%%%%%%%%%%%%%%%%%%%%%%%%%%%%%%%%%%
\section{Cotangent paths and their homotopy}%
\label{Cotangent-paths}                                 %
%%%%%%%%%%%%%%%%%%%%%%%%%%%%%%%%%%%%%%%%%%%%%%%%%%%%%%%%%%%%
%%%%%%%%%%%%%%%%%%%%%%%%%%%%%%%%%%%%%%%%%%%%%%%%%%%%%%%%%%%%

In the sequel $M$ will denote a Poisson manifold. We let $\#:T^*M\to
TM$ be the bundle map determined by contraction with the Poisson
tensor, and we let $[\cdot, \cdot]$ be the Lie the bracket on 1-forms,
which was defined in the introduction. We also denote by $\pi: T^*M\to
M$ the canonical projection. 

\begin{definition} A \textbf{cotangent path} in $M$ is a
path $a:I\to T^*M$, where $I=[0,1]$, satisfying
\[\frac{d}{dt}\pi(a(t))=\#a(t). \]
\end{definition}

Given a cotangent path $a:I\to T^*M$ and a vector field $X$
on $M$ we define, following \cite{GiGo}, the path integral:
\[ \int_a X\equiv\int_0^1 \langle a(t), X(\gamma(t))\rangle dt,\]
where $\gamma(t)=\pi(a(t)$ is the base path of $a$. If $X=X_h$ is a
Hamiltonian vector field for some Hamiltonian function $h\in
C^\infty(M)$, we see that
\[ \int_a X_h=h(\gamma(1))-h(\gamma(0)).\]
So for a Hamiltonian vector field the integral only depends on the
end-points of the base of the cotangent path. 

{From} a differential-geometric viewpoint, as was first explained in
\cite{Fer2}, cotangent paths are precisely the paths along which
parallel transport can be performed whenever a contravariant
connection has be chosen. Let us briefly recall how this works.

First of all, a \textbf{contravariant connection} $\nabla$ on a vector
bundle $E$ over $M$ can be thought of as a bilinear map
\[ \Omega^1(M)\times \Gamma(E)\to \Gamma(E),\ (\al,
s)\mapsto\nabla_{\al} s,\]  
satisfying the following identities:
\begin{enumerate}
\item[(a)] $\nabla_{f\al} s = f\nabla_{\al} s$,
\item[(b)] $\nabla_{\al}(fs)= f\nabla_{\al} s+\#\al(f) s$,
\end{enumerate}
for all $f\in C^{\infty}(M)$, $\al\in \Omega^1(M)$, $s\in
\Gamma(E)$. This is the analogue, in the Poisson category, of the usual
connections.

The value of $\nabla_{\al}s$ at $x\in M$ only depends on the
value of $\al$ at $x$ and the values of $s$ along the integral curve
of $\#\al$ through $x$. Also, just like the usual covariant derivative,
given a cotangent path $a$ and path $s$ in $E$ above $\pi\circ a$,
there is a well-defined contravariant derivative $\nabla_{a} s$, which 
is a path in $E$ above $\pi\circ a$. This leads to the
notion of \textbf{parallel transport} along a
cotangent path $a$: it is the map $\tau_a:E_{\pi(a(0))}\to
E_{\pi(a(1))}$ which takes $u\in E_{\pi(a(0))}$ to $s(1)\in
E_{\pi(a(1))}$, where $s:I\to E$ is the solution of the differential equation
$\nabla_{a} s = 0$, with initial condition $s(0)=u$. All this is
described in detail in \cite{Fer2}.

We will be mostly interested in contravariant connections on
$T^*M$. A covariant connection $\nabla$ on $TM$ induces, apart
from the dual covariant connection on $T^*M$, still denoted by
$\nabla$, a contravariant connection $\overline{\nabla}$ on $TM$
defined by
\[\overline{\nabla}_{\al}X=\nabla_{X}\#\al+ [\#\al,X],\]
and a dual contravariant connection on $T^*M$, also denoted by
$\overline{\nabla}$, and given by:
\[\overline{\nabla}_{\al}\be =\nabla_{\#\be}\al+[\al,\be].\]
These connections satisfy
$\overline{\nabla}\#=\#\overline{\nabla}$ and they form a
\textbf{basic connection} in the sense of \cite{Fer2}. The connection 
$\overline{\nabla}$ has \textbf{contravariant torsion} given by
\[ T_{\nabla}(\omega,\eta)= \overline{\nabla}_{\omega}\eta-
\overline{\nabla}_{\eta}\omega-[\omega,\eta].\]

Fix a covariant connection $\nabla$ on $TM$. Given a family
$a_{\epsilon}=a_{\epsilon}(t)$ of cotangent paths with the property
that the base paths $\gamma_{\epsilon}(t)=\pi(a_{\epsilon}(t))$
have fixed end points, we consider the differential
equation(\footnote{Given a covariant connection $\nabla$ on bundle
  $E\to M$, a path $\gamma(t)$ in $M$ and a path $u(t)$ in $E$ over
  $\gamma$, we use the notation 
$\partial_tu\equiv\nabla_\gamma u$.})
\begin{equation}
\label{eq-homotopy}
\partial_{t}b-\partial_{\epsilon}a=T_{\nabla}(a, b).
\end{equation}
This equation has a unique solution $b= b(\epsilon, t)$ (running
in $T^*M$) with initial condition $b(\epsilon,0)=0$, and we
define the \textbf{variation} of $a_{\epsilon}$ to be the cotangent path:
\[ \var(a_{\epsilon})= b(\epsilon, 1).\]

\begin{definition}
\label{homotopy} A \textbf{cotangent homotopy} is a family 
$a_{\epsilon}=a_{\epsilon}(t)$ of cotangent paths with the property 
that the base paths $\gamma_{\epsilon}(t)$ have fixed end points and
$\var(a_{\epsilon})= 0$. Two cotangent paths $a_0$ and $a_1$ are
said to be homotopic, and we write $a_0\sim a_1$ if there is a
cotangent homotopy joining them.
\end{definition}

For an alternative definition of cotangent homotopy which avoids the
use of connections we refer to \cite{CrFe}. 

Cotangent homotopy defines an equivalence relation on the set of
cotangent paths and enjoys 
 many properties similar to the ones
enjoyed by the usual notion of homotopy. For example, since the analogue of
closed 1-forms in this calculus are the Poisson vector fields, 
the following corresponds to the homotopy invariance of integrals of closed 1-forms
along paths:
\begin{proposition}
\label{hom-inv}
Let $a_0$ and $a_1$ be cotangent paths. If $a_0\sim a_1$ then, for
any Poisson vector field $X$ on $M$, we have
\[ \int_{a_0} X=\int_{a_1} X.\]
\end{proposition}

\begin{proof} 
Let $a=a(\epsilon,t)$ and $b=b(\epsilon,t)$ be as above, and
set $I=\langle a,X\rangle$, $J=\langle b,X\rangle$. Fixing some
connection $\nabla$, for any $\omega\in \Gamma(T^*M)$, we have
\[ \frac{d}{dt} \langle \omega(\gamma(t)), X(\gamma(t))\rangle =
\langle \nabla_{\frac{d\gamma}{dt}}\omega, X(\gamma(t))\rangle +
\langle\omega(\gamma(t)), \nabla_{\frac{d\gamma}{dt}}X\rangle .\]
A straightforward computation using the defining equation
(\ref{eq-homotopy}) and expression (\ref{kozul:bracket}) for the Lie
bracket, shows that
\[ \frac{dI}{d\epsilon}- \frac{dJ}{dt}= \Lie_X{\Pi}(a,b) .\]
Integrating first with respect to $t$, using $b(\epsilon, 0)= 0$, and
then integrating with respect to $\epsilon$, we find that
\[ \int_{a_{1}}X- \int_{a_{0}}X=
\int_{\var(a_\epsilon)} X + \int_{0}^{1}\int_{0}^{1}
\Lie_X \Pi (a, b)\, dtd\epsilon.\] for any vector field $X$, and any
family $\{a_{\epsilon}\}$.  For cotangent homotopies and Poisson vector
fields this gives $\int_{a_{0}}X= \int_{a_{1}}X$.
\end{proof}

\begin{remark}
In \cite{GiGo}, Ginzburg gives a simple example (cf.~Example 3.4)
showing that the integral is not invariant under homotopy of the base
path. Referring to the (non-)invariance of the integral under this
``naive homotopy'' he states: ``this example shows that the above
naive definition of homotopy is not a \emph{correct} extension of this
notion to the Poisson category''. The proposition above shows that
cotangent homotopy gives the desired extension.
\end{remark}

Another instance of behavior similar to the usual notion of homotopy
is obtained by looking at flat contravariant connections. For these
connections we have the following contravariant version of a classical
result in differential geometry:

\begin{proposition}
Given a flat contravariant connection, parallel transport
along cotangent paths is invariant under cotangent homotopy.
\end{proposition}

The proof is similar to the proof of the previous proposition (see
Proposition 1.6 in \cite{CrFe}), and will be omitted. Next we give the
contravariant analogue of the fundamental group.

\begin{definition}
Let $x\in M$ be a point in the Poisson manifold $M$. The
\textbf{isotropy group} $\Sigma(M, x)$ is the set of equivalence classes of
cotangent paths whose base paths in $M$ start and end at $x$. The
group structure is defined by concatenation of paths. We also
define the \textbf{restricted isotropy group} $\Sigma^0(M, x)$ by
considering only cotangent paths whose base path is a contractible loop.
\end{definition}

The groups $\Sigma(M, x)$ play a fundamental role in Poisson
geometry. For instance, one can show that linear holonomy of Poisson
manifolds (defined along cotangent paths, cf.~\cite{Fer2}), only
depends on homotopy classes, hence factors through the
  isotropy group.
Also, if $\mu: S\to M$ is a complete
symplectic realization, then there is a natural action of $\Sigma(M,
x)$ on $\mu^{-1}(x)$, and, if the ``reduction'' $\mu^{-1}(x)/\Sigma(M,
x)$ is smooth, then it carries a natural symplectic structure. More
details and more examples will be given in the later sections. Notice
also that, by Proposition \ref{hom-inv}, integration along closed
cotangent paths gives a group homomorphism
\[ \int:\Sigma(M, x)\to H^1_{\Pi}(M)^*,\]
where $H^1_{\Pi}(M)$ denotes the 1st Poisson cohomology group of
$M$, i.e., the set of Poisson vector fields modulo the Hamiltonian
vector fields. 

Though one may view the groups $\Sigma(M, x)$ as \emph{analogues} of the
fundamental groups of manifolds, there is also a strong analogy 
with the construction of simply connected Lie groups
integrating Lie algebras. In the remaining part of this section, we
elaborate on these two analogies. 

Notice that the two groups we have just defined do relate to the
fundamental groups of the symplectic leaf $L$ through $x$. We have
a short exact sequence
\begin{equation}
\label{eq1} \xymatrix{ 1\ar[r]& \Sigma^0(M, x)\ar[r]& \Sigma(M,
x)\ar[r]& \pi_{1}(L,x)\ar[r]& 1,}
\end{equation}
where the first map is just the inclusion, while the second one
associates to a cotangent path its underlying base path. In the
Poisson context however, there are several new features which are not
present in the case of classical fundamental groups.  For example, one
can assert that $\Sigma(M, x)$ and $\Sigma(M, y)$ are isomorphic only
when $x$ and $y$ lie in the same symplectic leaf of $M$. Also, the
groups $\Sigma(M, x)$ are not to be considered as discrete groups; 
this brings us to our second analogy, for which we need to recall the
reconstruction of simply connected Lie groups in terms of paths in
their Lie algebras, due to Duistermaat and Kolk \cite{DuKo}. 

Let $\widetilde{G}$ be a simply connected Lie group with Lie
algebra $\mathfrak{g}$. We can identify elements $g\in
\widetilde{G}$ with homotopy classes of paths $g(t)$ starting at
$g(0)=e$ and ending at $g(1)=g$. Using derivatives and right
translations, we obtain a bijective correspondence between paths
$g(t)$ in the Lie group $\widetilde{G}$ starting at $e$, and paths $a(t)$ in
the Lie algebra $\gg$ with free end points. Explicitly:
\begin{equation}
\label{derivatives}
a(t)=\left.\frac{d}{ds}g(s)g(t)^{-1}\right|_{s=t}.
\end{equation}
This bijection shows that there is a natural notion of homotopy of
paths in $\gg$, which corresponds to the homotopy of paths in
$\tilde{G}$ (with fixed end points): if $a_{\epsilon}$ is a family of
paths in $\gg$, the only thing we have to require is that the induced
family $g_{\epsilon}$ of paths in $\tilde{G}$ has fixed end points
$g_{\epsilon}(1)$. If we let $var(a_{\epsilon})=
\frac{d}{d\epsilon}g_{\epsilon}(1)$ one can show that
$var(a_{\epsilon})$ can be defined directly in terms of the family
$a_{\epsilon}$, i.e., using data in $\gg$ and with no reference to the
Lie group $\tilde{G}$. This is in complete analogy with the variation for
cotangent paths given above, with the extra simplification 
coming from
the fact that here no connection is involved: $var(a_{\epsilon})=
b(\epsilon, 1)$ where $b=b(\epsilon,t)$ is the solution of the
differential equation
\[ \frac{db}{dt}- \frac{da}{d\epsilon}=[b,a],\qquad b(\epsilon,0)=0.\]
The upshot is that the resulting group $G(\gg)$ of homotopy classes of
paths in $\gg$ can be defined directly in terms of the Lie algebra
$\gg$, and hence makes sense even if we don't assume the existence of
$\tilde{G}$. Moreover, since the space of paths in $\gg$ (say of class
$C^1$) carries a natural structure of Banach manifold, $G(\gg)$ has a
natural quotient topology, and at most one interesting smooth
structure (the one for which the quotient map is a submersion). The
fact that a smooth structure always exists, is far from trivial (see
\cite{DuKo}). This gives a proof of Lie's third theorem asserting the
existence of a simply-connected Lie group $G(\gg)$ integrating $\gg$,
and at the same time a proof of the fact that any simply-connected Lie
group $\tilde{G}$ must be isomorphic to $G(\gg)$.

Let us look now at the isotropy groups $\Sigma(M, x)$. Since the space
of cotangent paths carries a natural structure of a Banach manifold
(for which the underlying topology is the $C^1$-topology), $\Sigma(M,
x)$ has a natural quotient topology, and there is no ambiguity when
looking for the smooth structure on $\Sigma(M, x)$ (and similarly for
$\Sigma^0(M, x)$). However, unlike $G(\gg)$, the isotropy group
$\Sigma(M, x)$ is not always a Lie group (see Lemma \ref{lemma-Lie}
and Corollary \ref{mon-Lie}).

We recall from \cite{CrFe} why $\Sigma(M, x)$ may fail to be a Lie
group (see, also, the informal discussion in \cite{Sev}). Observe that
the Lie bracket on 1-forms induces a Lie bracket on the kernel of the
map $\#_x:T^*M\to TM$. By skew-symmetry, this kernel coincides with
the co-normal bundle to the leaf $L$ through $x$. We call
$\nu_{x}^{*}(L)=\Ker(\#_{x})$ the \textbf{isotropy Lie algebra} at
$x$. Now, the simply connected Lie group $G(\nu_{x}^{*}(L))$
integrating the isotropy Lie algebra and the restricted isotropy Lie
group $\Sigma^0(M, x)$ are both quotients of the space of paths in
$\nu_{x}^{*}(L)$. Also, two paths in $\nu_{x}^{*}(L)$ which are
homotopic as Lie algebra paths, are homotopic as cotangent
paths. However, the converse is not true, since a cotangent homotopy
may force one to go away from $x$. In other words, $\Sigma^{0}(M, x)$
is (topologically) a quotient of $G(\nu_{x}^{*}(L))$, and the
following result is immediate from \cite{CrFe}:

\begin{lemma}
\label{lemma-Lie} Let $M$ be  a Poisson manifold, $x\in M$, and
denote by $L$ the symplectic leaf through $x$. The following are
equivalent:
\begin{enumerate}
\item[(i)] $\Sigma(M, x)$ is Hausdorff;
\item[(ii)] $\Sigma(M, x)$ is a Lie group;
\item[(iii)] $\Sigma^0(M, x)$ is Hausdorff;
\item[(iv)] $\Sigma^0(M, x)$ is a Lie group.
\end{enumerate}
In this case, $\Sigma(M, x)$ has Lie algebra the isotropy Lie algebra
$\nu^{*}_{x}(L)$, its connected component of the identity is
$\Sigma^0(M, x)$, and the group of its connected components
$\pi_0(\Sigma(M, x))$ is isomorphic to $\pi_{1}(L, x)$.
\end{lemma}

One can get a better grasp on the notions of homotopy and variation
introduced above when $T^*M$ is integrable to a symplectic groupoid
$\Sigma$. Then, as in our discussion of Lie algebra paths, there is a
bijection between cotangent paths and paths $g(t)$ in $\Sigma$ which
stay in an $\s$-fiber, and which start at an identity element. The
variation of families $a_{\epsilon}$ of cotangent paths corresponds
then to the derivatives with respect to $\epsilon$ of the end points
of the associated paths in $\Sigma$. Hence, $a_{\epsilon}$ is a
cotangent homotopy precisely when the corresponding family of paths in
$\Sigma$ is a homotopy with fixed end points.

\begin{remark}
\label{remark-paths} 
Let us briefly explain the perspective one gains by using the language
of algebroids/groupoids, as well as its unifying role. First of all, a
cotangent path can be identified with a Lie algebroid morphism $TI\to
T^*M$. Similarly, two cotangent paths $a_0$ and $a_1$ are homotopic if
and only if there exists a Lie algebroid homomorphism $adt+bd\epsilon:
T(I\times I)\to T^*M$, which covers a (ordinary) homotopy with fix
end-points between the base paths $\pi(a_i(t))$, and which restricts
to $a_i(t)$ on the boundaries.  Therefore, all these constructions
readily extend to any Lie algebroid. For example, in the case $A=TM$
one recovers the usual notions of path and homotopy of paths, and the
fundamental group(oid). For a Lie algebra $\gg$, viewed as a Lie
algebroid over a point, one recovers the construction of Duistermaat
and Kolk of the simply connected Lie group $G(\gg)$\footnote{Of course, 
for Lie groups, one uses pointwise multiplication of paths instead of concatenation, 
as in any first course on fundamental groups.}

We point out that we have given above the precise meaning of
``cotangent homotopy with fix end-points''.  For example, a family of
cotangent paths $a_\eps(t)$, $\eps\in[0,1]$, such that
$a_\eps(0)=a_\eps(1)=0$, in general, will \emph{not be} a cotagent
homotopy. In fact, such a definition would quickly lead to erroneous
statements: for example, integration of Poisson vector fields over
cotagent paths would not be invariant under cotangent
homotopy. Similar remarks apply for Lie algebroids. We refer the
reader to \cite{CrFe} for details on these constructions in the
general Lie algebroid framework.
\end{remark}

%%%%%%%%%%%%%%%%%%%%%%%%%%%%%%%%%%%%%%%%%%%%%%%%%%%
%%%%%%%%%%%%%%%%%%%%%%%%%%%%%%%%%%%%%%%%%%%%%%%%%%%
\section{The monodromy groups of Poisson manifolds}                    %
\label{monodromy}                                 %
%%%%%%%%%%%%%%%%%%%%%%%%%%%%%%%%%%%%%%%%%%%%%%%%%%%
%%%%%%%%%%%%%%%%%%%%%%%%%%%%%%%%%%%%%%%%%%%%%%%%%%%

In this section we give various descriptions of the monodromy
groups of a Poisson manifold $M$ at a point $x\in M$. We shall see
that these are certain additive subgroups
\[ \NN_{x} \subset \nu^*(L)_x ,\]
where $L$ is the symplectic leaf through $x$, and $\nu_{x}^{*}(L)$
the co-normal space to $L$ at the point $x$. Recall that
$\nu_{x}^{*}(L)$ carries a Lie algebra structure induced from the
Lie bracket on 1-forms. The monodromy group actually sits inside
the center of this Lie algebra.

We start with the shortest possible description of monodromy:

\begin{descr}
The monodromy group $\NN_{x}$ is the set of vectors $v\in
Z(\nu_{x}^{*}(L))$, with the property that the constant cotangent
path $a(t)= v$ is cotangent homotopic to the zero cotangent path.
\end{descr}

To understand why this group shows up when looking at smoothness
issues, we consider a slightly larger group: $\tilde{\NN}_{x}\subset
G(\nu_{x}^{*}(L))$ consists of those elements represented by Lie
algebra paths which, as cotangent paths, are homotopic to the zero
path. First of all, by its very definition, there exists an exact sequence
\begin{equation}
\label{eq2} \xymatrix{ 0\ar[r]& \tilde{\NN}_{x}\ar[r]&
  G(\nu_{x}^{*}(L)) \ar[r]& \Sigma^0(M, x) \ar[r]& 0,}
\end{equation}
so that:
\[ \Sigma^0(M, x)= G(\nu_{x}^{*}(L))/\tilde{\NN}_{x}.\] 
Next, the group $\tilde{\NN}_{x}$ almost coincides with $\NN_{x}$. The
only difference comes from the fact that, although $\tilde{\NN}_{x}$
does lie in the center $Z(G(\nu_{x}^{*}(L)))$, it is only the
connected component $Z^{0}(G(\nu_{x}^{*}(L)))$ of the center that
naturally identifies (via the exponential) with the center of
$\nu_{x}^{*}(L)$. It follows that:
\[ \NN_{x}= \tilde{\NN}_{x}\cap Z^{0}(G(\nu_{x}^{*}(L))) \]

In particular, let us point out the following (see also Lemma
\ref{lemma-Lie})

\begin{corollary}
\label{mon-Lie} The isotropy group $\Sigma(M, x)$ is a Lie group
if and only if the monodromy group $\NN_{x}$ is discrete.
\end{corollary}

The next description not only explains the nature of the
monodromy, but it is also suitable for computations.

\begin{descr}
\label{descr2}
Consider a linear splitting $\sigma: TL\to T^*_{L}M$ of the short
exact sequence
\begin{equation}\label{seseq}
\xymatrix{0\ar[r]& \nu^*(L)\ar[r]& T^*_{L}M\ar[r]^{\#}& TL\ar[r]&
0}
\end{equation}
Assume that $\sigma$ can be chosen so that its curvature
\[ \Omega_{\sigma}(X, Y)\equiv \sigma ([X, Y])- [\sigma(X), \sigma(Y)]\]
is $Z(\nu^*(L))$-valued. Then
\[ \NN_{x}= \set{ \int_{\gamma} \Omega_{\sigma}:
[\gamma]\in \pi_{2}(L, x)}\subset \nu^*(L)_x.\]
\end{descr}
\vskip 10 pt

The integrals involved in this last description should be viewed as
integrals of forms with coefficients in flat vector
bundles(\footnote{Note that forms with values in flat 
  bundles can be integrated over simply connected cycles.}): the Bott
connection gives $\nu(L)$ (hence also $\nu^*(L)$, and $Z(\nu^*(L))$) a
natural flat vector bundle structure over $L$. We point out that the
assumption that $\sigma$ is center-valued is only needed to simplify
the outcome. It is satisfied in many examples (e.g.~whenever $x$ is a
regular point), but not always (see Example \ref{su(3)} in Section
\ref{Integrability:general}). In general, the formula
\begin{equation}
\label{eq:boundary}
\partial([\gamma])=
\int_0^1
\Omega_{\sigma}\left(\frac{d\gamma}{dt},\frac{d\gamma}{d\epsilon}\right)
d\epsilon
\end{equation}
defines a group homomorphism 
\[ \partial:\pi_{2}(L,x)\to G(\nu^*(L)_x)\ ,\]
and
\[ \tilde{\NN}_{x}= \{ \partial(\gamma) :
[\gamma]\in \pi_{2}(L, x)\} \subset G(\nu^*(L)_x).\] 
In the center-valued case, any center-valued Lie algebra path can
be identified (i.e. defines the same Lie group element) with its integral, viewed as a
central Lie group
element, and we obtain
\[ \partial([\gamma])=
\int_0^1\int_0^1
\Omega_{\sigma}\left(\frac{d\gamma}{dt},\frac{d\gamma}{d\epsilon}\right)
d\epsilon dt=\int_{\gamma} \Omega_{\sigma}\ ,
\]
which is the description of the $\NN_{x}$ given above.

Let us give a more conceptual explanation for the homomorphism
$\partial$, which also provides the relation between the previous two
descriptions.  The ``philosophical idea'' on Lie algebroids mentioned
in the introduction, and the remark that the homotopy long exact
sequence of a (smooth) fibration can be carried out at the level of
tangent bundles (see the construction of $\partial$ below), suggests a construction of such sequences for
algebroids.

\begin{descr}
Viewing (\ref{seseq}) as analogous
to a fibration, there is a long exact sequence
\[ \xymatrix{\cdots\ar[r]&\pi_{2}(L,x)\ar[r]^{\partial}& G(\nu^*(L)_x)
\ar[r]^{j}& \Sigma(M, x)\ar[r]& \pi_1(L,x) },\] 
where $\partial:\pi_{2}(L,x)\to G(\nu^*(L)_x)$ is the map given above, 
and $j$ is the composition of the projection from $G(\nu^*(L)_x)$ onto
$\Sigma^{0}(M, x)$ with the inclusion.
\end{descr}

One can construct $\partial$ in complete analogy with the
construction of the boundary map for fibrations: given a $2$-loop
$\gamma: I\times I\to L$ based at $x$, one lifts its differential
$d\gamma$ to a morphism $adt+bd\epsilon:T(I\times I)\to A$. Since
$\gamma$ on the boundary of $I\times I$ takes the constant value
$x$, the restriction of this morphism to the boundary, gives the
path $a(1,t)$ in the Lie algebra $\nu^*(L)_x$, hence an element in
the Lie group $G(\nu^*(L)_x)$. If we choose a splitting $\sigma$,
we can use it to lift $d\gamma$, and the computation gives
precisely the expression (\ref{eq:boundary}) for $\partial$.

%%%%%%%%%%%%%%%%%%%%%%%%%%%%%%%%%%%%%%%%%%%%%%%%%%%
%%%%%%%%%%%%%%%%%%%%%%%%%%%%%%%%%%%%%%%%%%%%%%%%%%%
\section{The symplectic groupoid $\Sigma(M)$}%
\label{symplectic-groupoid}
%%%%%%%%%%%%%%%%%%%%%%%%%%%%%%%%%%%%%%%%%%%%%%%%%%%
%%%%%%%%%%%%%%%%%%%%%%%%%%%%%%%%%%%%%%%%%%%%%%%%%%%

In this section we describe the Weinstein groupoid $\Sigma(M)$
associated to a Poisson manifold. Note that, though $\Sigma(M)$ is
just a special case (with $A=T^*M$) of the general construction of the
Weinstein groupoid of a Lie algebroid $A$ \cite{CrFe}, the
extra structure coming from the Poisson bracket and the cotangent
bundle translates into some special properties of $\Sigma(M)$. In
particular, we explain the relation with the work of Cattaneo and
Felder \cite{CaFe} on the Poisson-sigma model (the Hamiltonian
formulation).

We denote by $P(T^*M)$ the space of cotangent paths in $M$ which are
of class $C^1$, and we denote by $\sim$ the cotangent homotopy
equivalence defined using cotangent homotopies $a_{\epsilon}$ in
$P(T^*M)$ which are of class $C^2$ in $\epsilon$. The Weinstein
groupoid is defined as the quotient(\footnote{In \cite{CrFe}, this
groupoid was denoted $\G(T^*M)$.})
\begin{equation}
\label{eq:weinstein:groupoid} \Sigma(M)\simeq P(T^*M)/\sim\ .
\end{equation}
For the groupoid structure, the source and target maps
$\s,\t:\Sigma(M)\to M$ take the equivalence class of a cotangent path
to the end-points of the base path, while multiplication is given by
concatenation (see \cite{CrFe} for details). We
shall see that this groupoid gives new insights into the global
properties of a Poisson manifold.

As a simple example, let us consider the integration of Poisson
vector fields along cotangent paths, discussed in Section
\ref{Cotangent-paths} above.  First, invariance under cotangent
homotopy (see Proposition \ref{hom-inv}) shows that we can view
integration as a map which associates to each Poisson vector field $X$
a map
\[ \int X: \Sigma(M)\to \mathbb{R},\qquad  g=[a]\mapsto \int_{a}X.\]
Second, additivity with respect
to the concatenation of paths
\[ \int_{g_{0}g_{1}}X= \int_{g_{0}}X+\int_{g_{1}}X,\] 
can be viewed as a cocycle condition for $\int X$: if $\Sigma(M)$ is a
smooth manifold, this means that $\int X$ defines a \textbf{differentiable
1-cocycle} on $\Sigma(M)$. Third, and last, for an Hamiltonian
vector field we can write:
\[ \int_{g}X_{h}= h(\t(g))- h(\s(g)),\]
which just says that $\int X_h$ is a \textbf{differentiable
  coboundary}. Hence, denoting by $H^{1}_{\text{dif}}(\Sigma(M))$ the first
  differentiable cohomology group of $\Sigma(M)$, integration defines
  a map $H^{1}_{\Pi}(M)\to H^{1}_{\text{dif}}(\Sigma(M))$.

On the other hand, there is a Van Est map (see \cite{Cra}) which takes a
differentiable cocycle into a Poisson cohomology class, and which is
defined in all degrees. This map is known to be an isomorphism in
degree one. One can think of the Van Est map as differentiation of
differentiable cocycles. A simple check shows that:

\begin{proposition}
The map $\int: H^{1}_{\Pi}(M)\to H^{1}_{\text{dif}}(\Sigma(M))$ is the inverse
of the Van Est map in degree one. 
\end{proposition}

The reader will notice that all this discussion extends to any Lie
algebroid with appropriate changes.

On $\Sigma(M)$ we consider always the quotient topology, so that it
becomes a topological groupoid. As in the case of the isotropy groups
$\Sigma(M,x)$, there is no ambiguity when looking at the smoothness of
$\Sigma(M)$: the smooth structure, if exists, will be the unique one
for which the quotient map $P(T^*M)\to \Sigma(M)$ is a submersion. It
is also convenient to consider the larger space $\widetilde{P}(T^*M)$
of all $C^1$-curves $a: I\to T^* M$, with base path $\gamma= \pi\circ
a$ of class $C^2$. It has an obvious structure of a Banach manifold,
and $P(T^*M)$ is a (Banach) submanifold of $\widetilde{P}(T^*M)$
(cf. Lemma 4.6 in \cite{CrFe}). We need an explicit description of the
tangent spaces to these Banach manifolds, and a more geometric
description of the equivalence relation defined on them by cotangent
homotopies.

The tangent space $T_{a}\widetilde{P}(T^*M)$ consists of curves
$U:I\to TT^*M$ such that $U(t)\in T_{a(t)}T^*M$, i.e.. vector
  fields along $a$. Using a connection
$\nabla$ on $TM$ and the associated contravariant connections
$\overline{\nabla}$ (see Section \ref{Cotangent-paths}), such a curve
can be viewed as a pair $(u,\phi)$ consisting of curves $u:I\to T^*M$
and $\phi:I\to TM$ over $\gamma$ (namely, the vertical and horizontal
component of $U$). The subspace $T_{a}P(T^*M)\subset
T_{a}\widetilde{P}(T^*M)$ consists (see \cite{CrFe}, Section 5.2) of
those pairs $U= (u,\phi)$ with the property that
\[
\#u= \overline{\nabla}_{a}\phi.
\]

Now let us describe the equivalence relation defined by cotangent
homotopies in terms of a Lie algebra action. The Lie algebra is formed
by time-dependent 1-forms, vanishing at the end-points
\[
P_{0}\Omega^1(M)=\set{\eta_{t}\in \Omega^1(M), t\in I:\eta_{0}
                =\eta_{1}=0,\ \eta_t\text{ of class }C^1\text{ in }t}
\]
with Lie bracket the bracket on 1-forms with time varying as a
parameter. The Lie algebra $P_{0}\Omega^1(M)$
acts on $\widetilde{P}(T^*M)$ and this action satisfies the following
properties: 
\begin{enumerate}
\item[(a)] it is tangent to $P(T^*M)$; 
\item[(b)] the orbits in $P(T^*M)$ define a finite codimension
  foliation;
\item[(c)] two cotangent paths are homotopic if and only if they
belong to the same orbit. 
\end{enumerate}
This is proved in \cite{CrFe}, and here we shall give only the
definition of the action. This means we will define a Lie algebra map
\[  P_{0}\Omega^1(M)\to\X(\widetilde{P}(T^*M)),\quad \eta\mapsto X_{\eta} \]
So, given a time-dependent 1-form $\eta\in P_{0}\Omega^1(M)$ and a
path $a\in \widetilde{P}(T^*M)$ with underlying path $\gamma:I\to M$, we
have to describe a path $X_{\eta,a}$ in $T^*M$ above $\gamma$.
Let us first assume that $a$ is a cotangent path. To describe
$X_{\eta, a}$, we specify the components $(u, \phi)$ with respect
to a connection $\nabla$:
\[
u= \overline{\nabla}_{a}b,\quad \phi= \#b.
\]
This does not depend on the connection, and, by the discussion
above, it does define a vector tangent to $P(T^*M)$. We check that
this can
also be written as:
\[
X_{\eta, a}(t)= {\left. \frac{d}{d\epsilon}\right|}_{\epsilon= 0}
\phi_{\eta}^{\epsilon, 0} a(t)+ \frac{d\eta_{t}}{dt} (\gamma(t))
\]
where $\phi_{\eta}^{s, t}$ stands for the flow of the time
dependent one form $\eta$ (\footnote{Given a time-dependent one-form
$\eta$ on $M$, its flow is a map $\phi_{\eta}^{t, s}(x):
T^{*}_{x}(M)\to T^{*}_{\Phi_{\#\al}^{t,s}}(M)$. It is uniquely determined
by the conditions
$\phi^{t,s}_{\eta}\phi^{s,u}_{\eta}=\phi^{t,u}_{\eta}$, 
$\phi^{t,t}_{\eta}=$id, and the differential equation 
\[{\left. \frac{d}{dt}\right|}_{t=s} (\phi_{\eta}^{t,s})^{*}\beta
=[\eta^s,\beta],\]
where $(\phi_{\eta}^{t, s})^{*}(\beta)(x)=
\phi_{\eta}^{s,t}\beta(\phi_{\#\eta}^{t, s}(x))$, with
$\phi_{\#\eta}^{t,s}$ the flow of the vector field $\#\eta$.}).  
Now, this formula make sense for \emph{any} path $a:I\to T^* M$
and this is our definition of the action.

What we have said so far makes sense for any Lie algebroid. Let us now
take advantage of the fact that we have $A=T^*M$, the cotangent
bundle. We have a natural identification $\widetilde{P}(T^*M)\simeq
T^*P(M)$, where $P(M)$ denotes the Banach space of paths $\gamma:I\to
M$ of class $C^2$. Hence, $\widetilde{P}(T^*M)$ carries a natural
symplectic structure. Moreover, the set of cotangent paths $P(T^*M)$
is the level set $J^{-1}(0)$ of the map
$J:\widetilde{P}(T^*M)\to P_{0}\Omega^1(M)^*$ given by:
\[ \langle J(a),\eta\rangle=
\int_0^1 \langle
\frac{d}{dt}\pi(a(t))-\#a(t),\eta(t,\gamma(t))\rangle dt.\]
and one has the following result due to Cattaneo and Felder (see
\cite{CaFe}):

\begin{theorem}
The Lie algebra action of $P_{0}\Omega^1(M)$ on
$\widetilde{P}(T^*M)$ is Hamiltonian, with equivariant moment map
$J:\widetilde{P}(T^*M)\to P_{0}\Omega^1(M)^*$.
\end{theorem}

Hence the Weinstein groupoid can be described alternatively as a
Marsden-Weinstein reduction:
\begin{equation}
\label{eq:sigma:model}
\Sigma(M)=\widetilde{P}(T^*M)//P_{0}\Omega(M).
\end{equation}

The two alternate descriptions (\ref{eq:weinstein:groupoid}) and
(\ref{eq:sigma:model}) for the Weinstein groupoid of a Poisson
manifold give the precise relationship between our integrability
approach introduced in \cite{CrFe}, which is valid for any Lie
algebroid, and the approach of Cattaneo and Felder in \cite{CaFe}.

\begin{corollary}
If $\Sigma(M)$ is smooth, then it admits a symplectic form which
turns $\Sigma(M)$ into a symplectic groupoid.
\end{corollary}

\begin{proof}
We only need to check the compatibility of the symplectic form
with the product, i.e., that the graph of multiplication is a Lagrangian
submanifold. This can be restated in a slightly simpler form using the 
following proposition which will also be useful later.

\begin{proposition}
\label{prop:compatible:symplectic}
Let $\G$ be a Lie groupoid, and let $\omega\in\Omega^2(\G)$ be a
symplectic form. The following statements are equivalent:
\begin{enumerate}
\item[(i)] The graph of multiplication $\gamma_m=\{(g,h,g\dot h)\in
    \G\times\G\times\G:(g,h)\in\G^{(2)}\}$ is a Lagrangian submanifold
  of $\G\times\G\times\bar{\G}$;
\item[(ii)] The relation $m^*\omega=\pi^*_1\omega+\pi^*_2\omega$ holds;
\end{enumerate}
where we denote by $m:\G^{(2)}\to\G$ the multiplication in $\G$ and by
$\pi_1,\pi_2:\G^{(2)}\to\G$ the projections to the first and second factors.
\end{proposition}

Using this Proposition we prove the Corollary above. First note that
we have the following explicit formula for the symplectic form
$\widetilde{\omega}$ in $\widetilde{P}(T^*M)$:
\[ \widetilde{\omega}_a(U_1, U_2)=
\int_0^1 \omega_{can}(U_1(t), U_2(t))dt,\]
for all  $U_1, U_2\in T_a\widetilde{P}(T^*M)$, where $\omega_{can}$ 
is the canonical symplectic form on $T^*M$. The additivity of the integral
shows that that condition (ii) holds at the level of
$\widetilde{P}(T^*M)$, hence it must hold also on the reduced space
$\Sigma(M)$.
\end{proof}

\begin{proof}[Proof of Proposition \ref{prop:compatible:symplectic}]
We claim that the graph of multiplication $\gamma_m$ satisfies the
following two properties 
\begin{enumerate}
\item[(a)] $\gamma_m$ is isotropic iff
  $m^*\omega=\pi^*_1\omega+\pi^*_2\omega$;
\item[(b)] $\gamma_m$ is isotropic implies that $M$ is also isotropic;
\end{enumerate}
Assuming these claims, it should be clear that (i) implies
(ii). Conversely, if (ii) holds, then by (a) we have that $\gamma_m$
is an isotropic submanifold. Since
$\dim\gamma_m=\dim\G^{(2)}=2\dim\G-\dim M$, we see that $\gamma_m$ is
Lagrangian provided that $\dim M=\frac{1}{2}\dim\G$.  Since $\gamma_m$ is 
isotropic, we know already that $2\dim\G-\dim M\le \frac{3}{2}\dim\G$,
and so we have
\[ \dim M\ge \frac{1}{2}\dim\G.\]
But, by (b), $M$ is also an isotropic submanifold
of $\G$ so that
\[ \dim M\le \frac{1}{2}\dim\G.\]
and equality must indeed hold. 

To prove the claims let us  denote by
$\Omega=\omega\oplus\omega\oplus(-\omega)$ the symplectic 
form on $\G\times\G\times\bar{\G}$, and let
$\gamma:\G^{(2)}\to\G\times\G\times\G$ be the embedding:
\[(g,h)\mapsto (g,h,g\cdot h).\]
Obviously, the graph $\gamma_m$ is an isotropic submanifold of
$\G\times\G\times\bar{\G}$ if and only if $\gamma^*\Omega=0$, and this condition
is clearly equivalent to  
\[ m^*\omega=\pi^*_1\omega+\pi^*_2\omega.\]
This shows that (a) holds.

In order to show that (b) also holds we observe that $M$, viewed as
the identity section of $\G$, can be expressed as
\[ M=\set{x\in\G: \exists g\in \G, x\cdot g=g}.\]
Hence we have $M=\pi(\Delta\cap\gamma_m)$ where
$\pi:\G\times\G\times\G\to\G$ is projection onto the first fact and 
\[\Delta\equiv\set{(h,g,g)\in\G\times\G\times\G}.\]
Now, it is easy to check that $\Delta\subset \G\times\G\times\bar{\G}$ is a
coisotropic submanifold, $\pi:\Delta\to\G$ is projection along its
characteristic foliation, and $\Delta$ and $\gamma_m$ have clean
intersection. By a standard argument in symplectic geometry (see
\cite{MaSa}, Lemma 5.34), $\gamma_m$ being isotropic implies that the
projection $M=\pi(\Delta\cap\gamma_m)$ is also isotropic.
\end{proof}

%%%%%%%%%%%%%%%%%%%%%%%%%%%%%%%%%%%%%%%%%%%%%%%%%%%
%%%%%%%%%%%%%%%%%%%%%%%%%%%%%%%%%%%%%%%%%%%%%%%%%%%
\section{The smoothness of $\Sigma(M)$}%
\label{smoothness}
%%%%%%%%%%%%%%%%%%%%%%%%%%%%%%%%%%%%%%%%%%%%%%%%%%%
%%%%%%%%%%%%%%%%%%%%%%%%%%%%%%%%%%%%%%%%%%%%%%%%%%%

We are now in position to give a proof of our main result concerning
the smoothness of $\Sigma(M)$. The discussion here depends heavily on
our general integrability criteria for Lie algebroids \cite{CrFe}.

\begin{theorem}
\label{strong-form} 
For a Poisson manifold $M$ the following statements are equivalent: 
\begin{enumerate}
\item[(i)] $M$ is integrable by a symplectic Lie groupoid,
\item[(ii)] the algebroid $T^*M$ is integrable,
\item[(iii)] the Weinstein groupoid $\Sigma(M)$ is a smooth manifold, 
\item[(iv)] the exponential map of $\Sigma(M)$ associated to a (or any)
connection is injective around the zero section, 
\item[(v)] the monodromy groups $\NN_{x}$, where $x\in M$, are locally  uniformly
discrete.
\end{enumerate}
In this case $\Sigma(M)$ is the unique $s$-simply connected,
symplectic groupoid which integrates $M$.
\end{theorem}

Before proceeding to the proof, we clarify conditions (iv) and (v).

The exponential map referred to in condition (iv) is defined as
follows. Given a contravariant connection $\nabla$, one can define the
notion of cotangent geodesic as in the classical case, by the equation
$\nabla_{a}a= 0$ (see \cite{Fer2}). Then, for $a_{0}\in
T^*M$ sufficiently close to zero, the geodesic starting at $a_0$ is
defined up to the time $t=1$, hence defines an element in
$P(T^*M)$. This correspondence defines a map
$\widetilde{\exp}_{\nabla}: T^*M\to P(T^*M)$, and, the map induced to 
the quotient, $\exp_{\nabla}: T^*M\to \Sigma(M)$, is called the
\textbf{exponential map} associated to $\nabla$. As in the classical
case, the map is only defined on a neighborhood of the zero section.

Condition (v) should be viewed as the main ``integrability
criterion''. 
The results of Section \ref{monodromy} show that it is suitable for
computation in concrete examples. In order to measure the discreteness of the
monodromy groups we define a function $r_{\NN}:M\to [0,+\infty]$ by
\[ r_{\NN}(y)= d(\NN_{y}-\{0_{y}\}, \{0_{y}\}),\]
where the distance is induced by some norm on $T^*M$. By convention,
$r_{\NN}(y)=+\infty$ if $\NN_{y}=\{0_{y}\}$. Clearly, $\NN_{y}$ is
discrete if and only if $r_{\NN}(y)\neq 0$. By uniform discreteness of
$\NN$ around $y$ we mean that
\[ \liminf_{y\to x}r_{\NN}(y)\neq 0.\]

\begin{proof}[Proof of Theorem \ref{strong-form}]

If $\Sigma$ is a symplectic groupoid integrating $M$, then its
associated algebroid is isomorphic to $T^*M$, and this shows that (i)
implies (ii). To see that (ii) implies (iii), we assume that $\Sigma$
is a Lie groupoid integrating $T^*M$.  Taking the universal covers of
the $\s$-fibers of $\Sigma$ together, we get a new groupoid
$\widetilde{\Sigma}$ which still integrates $T^*M$, and which is
$\s$-simply connected. More precisely, $\widetilde{\Sigma}$ is the
collection of paths lying on the $\s$-fibers, which start at the
identity section, where we identify any two which are homotopic by
homotopies with fixed end points (see \cite{CrFe} for details). It
follows from Remark \ref{remark-paths}
 that $\tilde{\Sigma}$ is
naturally identified with the space of cotangent paths, modulo
homotopy, i.e. with $\Sigma(M)$. Since the quotient $\Sigma$ of
$\tilde{\Sigma}$ is smooth, so is $\tilde{\Sigma}\cong
\Sigma(M)$. That (iii) implies (i) is the content of the previous
section.

For any contravariant connection $\nabla$ the exponential
$\widetilde{\exp}_{\nabla}: T^*M\to P(T^*M)$ 
defines a
transversal to the foliation on $P^*(TM)$. In particular, if the
associated leaf space $\Sigma(M)$ is smooth, the differential of
$\exp_{\nabla}$ at zero elements is the invertible, and this shows that
(iii) implies (iv). Conversely, if (iv) holds, then one uses
$\exp_{\nabla}$ to define the smooth structure on $\Sigma(M)$ around
the identity elements, and then, using right translations, one extends
the smooth structure to the entire topological groupoid
$\Sigma(M)$ (cf. \cite{CrFe}). This explains why (iv) implies (iii).

Finally, the equivalence with (v) is the main theorem in \cite{CrFe} ,
applied to the special case
of the cotangent Lie algebroid of a Poisson
manifold.
\end{proof}

%%%%%%%%%%%%%%%%%%%%%%%%%%%%%%%%%%%%%%%%%%%%%%%%%%%
%%%%%%%%%%%%%%%%%%%%%%%%%%%%%%%%%%%%%%%%%%%%%%%%%%%
\section{The regular case}%
\label{Integrability:regular}                     %
%%%%%%%%%%%%%%%%%%%%%%%%%%%%%%%%%%%%%%%%%%%%%%%%%%%
%%%%%%%%%%%%%%%%%%%%%%%%%%%%%%%%%%%%%%%%%%%%%%%%%%%

A Poisson manifold is called regular if the Poisson tensor has
constant rank. In the regular case, both $\F=\Im\#$ and
$\nu^*(\F)=\ker\#$ are vector bundles over $M$, and the isotropy
Lie algebras $\nu^*(\F)_{x}$ are all abelian. This amounts to
several simplifications and leads to a new geometric
interpretation of the integrability criteria.

Let us consider, for instance, the long exact sequence of Section
\ref{monodromy}. In the regular case it becomes
\[ \xymatrix{\cdots\ar[r]&\pi_{2}(L,x)\ar[r]^{\partial}& \nu^*(L)_x
\ar[r]& \Sigma(M, x)\ar[r]& \pi_1(L,x) }.\] The monodromy groups
can be described (or defined) as the image of $\partial$, which,
in turn, is given by integration of a canonical cohomology 2-class
$[\Omega_{\sigma}]\in H^2(L,\nu_{x}^{*}(L))$, and this class can
be computed explicitly by using a section $\sigma$ of $\#:
T^{*}_{L}M\to TL$. In this section we will give a different
geometric interpretation of monodromy based on the symplectic
geometry of the leaves of $M$.  We state the main results first,
and delay all proofs until
the end of the section.

Fix a point $x$ in a Poisson manifold $M$, let $L$ be the symplectic
leaf through $x$, and consider a 2-sphere $\gamma: S^2\to L$, which
maps the north pole $p_{N}$ to $x$. The symplectic area of $\gamma$ is
given, as usual, by
\[ A_{\omega}(\gamma)= \int_{S^2} \gamma^{*}\omega ,\]
where $\omega$ is the symplectic 2-form on the leaf $L$. By a
\textbf{deformation} of $\gamma$ we mean a family $\gamma_{t}: S^2\to M$
of 2-spheres parameterized by $t\in (-\eps,\eps)$, starting at
$\gamma_{0}= \gamma$, and such that for each fixed $t$ the sphere
$\gamma_t$ has image lying entirely in a symplectic leaf. The
\textbf{transversal variation} of $\gamma_t$ (at $t=0$) is the class of
the tangent vector
\[
\var_{\nu}(\gamma_{t})\equiv
\left[\left.\frac{d}{dt}\gamma_{t}(p_{N})\right|_{t=0}\right]\in
\nu_{x}(\F).
\]

\begin{figure}[h]
\begin{center}
        \epsfxsize=\textwidth
        \leavevmode
        \epsfbox{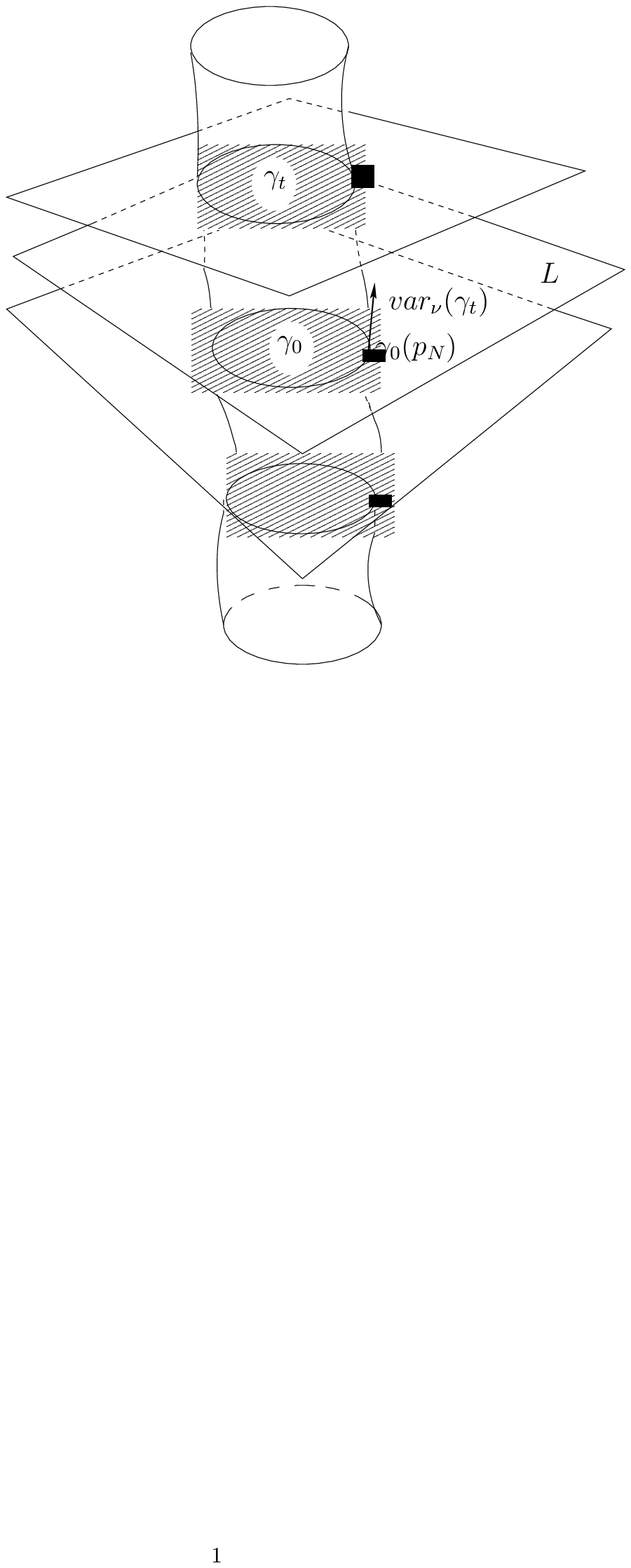}
\end{center}
        \caption{Symplectic spheres.}
\end{figure}

We shall see below that the quantity
\[ \left.\frac{d}{dt}A_{\omega}(\gamma_{t})\right|_{t=0} \]
only depends on the homotopy class of $\gamma$ and on
$\var_{\nu}(\gamma_{t})$. Finally, the formula
\[ \langle A^{'}_{\omega}(\gamma), var_{\nu}(\gamma_{t})\rangle=
\left.\frac{d}{dt}A_{\omega}(\gamma_{t})\right|_{t=0},\] 
applied to different deformations of $\gamma$, gives a well defined
element
\[ A^{'}_{\omega}(\gamma)\in \nu_{x}^{*}(L) .\]
The promised result is:

\begin{proposition}\label{var-area}
For any regular manifold $M$, the function $r_\NN$ is lower
semi-continuous, and
\[ \NN_{x}= \{A^{'}_{\omega}(\sigma): \sigma\in \pi_{2}(L, x)\}.\]
\end{proposition}

Next, we form the (set theoretical) bundle of ``variations of
symplectic areas'' (this already appears in \cite{AlHe} and in earlier work of Dazord, and 
it has been used in \cite{Xu2} in the case where the symplectic foliation is simple):
\[ \AAA^{'}(M)= \bigcup \NN_{x} \subset \nu^{*}(\F) .\]
Similarly, the groups $\S_{x}(M)= \nu_{x}^{*}(L)/\NN_{x}$ fit
together into a (set theoretical, again) bundle of groups,
\[ 
\S(M) =  \nu^*(\F)/\AAA^{'}(M)= \bigcup_{x\in M} \nu^*(\F)_{x}/ \NN_{x},
\]
which will be called the \textbf{structure groupoid} of $M$. Note
that the groups $\S_{x}(M)$ are closely related to the isotropy
groups $\Sigma(M, x)$ introduced in Section \ref{Cotangent-paths}.
More precisely, it follows from Section \ref{monodromy} that:
\[
\S_{x}(M)= \Sigma^{0}(M, x),
\] 
and we have a short exact sequence
\[
\xymatrix{1\ar[r]& \S_{x}(M)\ar[r]& \Sigma(M, x)\ar[r]&
  \pi_1(L)\ar[r]& 1.}
\]

In the sequel we will only be interested on the structures on
$\S(M)$ coming from $\nu^{*}(\F)$ via the projection
$\nu^{*}(\F)\to \S(M)$, and on structures on $\AAA^{'}(M)$ as a
subspace of $\nu^{*}(\F)$. In particular, both $\AAA^{'}(M)$ and
$\S(M)$ are topological groupoids, and there is no ambiguity when
asking for \emph{the} smooth structure on $\AAA^{'}(M)$ and on
$\S(M)$, respectively.

\begin{theorem}
\label{structure} 
For a regular Poisson manifold $M$, the following are equivalent:
\begin{enumerate}
\item[(i)] $M$ is integrable; 
\item[(ii)] $\Sigma(M)$ is a {\rm Lie} groupoid; 
\item[(iii)] $\S(M)$ is a {\rm Lie} groupoid;
\item[(iv)] $\AAA^{'}(M)$ is a {\rm Lie} groupoid.
\end{enumerate}
Moreover, in this case, $\AAA^{'}(M)$ is \'etale (i.e.~its source
is a local diffeomorphism), and there is an exact sequence of Lie
groupoids
\[ 
\xymatrix{ 1\ar[r]& \S(M)\ar[r]& \Sigma(M)\ar[r]^{\Phi}&
  \pi_1(\F)\ar[r]&1}.
\]
\end{theorem}

Here $\pi_1(\F)$ stands for the \textbf{homotopy groupoid} of $\F$.
This object is well known in foliation theory (sometimes under the
name of monodromy groupoid), and it is one of the simplest
examples of a Weinstein groupoid (namely the one associated to the
Lie algebroid $T\F$). The arrows between $x, y\in M$ are the
leafwise homotopy classes of paths (paths tangent to $\F$) with
end-points in $x$ and $y$, so that there are arrows only between
points lying in the same leaf.

In the short exact sequence above, $\S(M)$ and $\pi_{1}(\F)$ are
(in principle) easily computable from the Poisson geometry of $M$,
and the sequence will give in fact a pretty good indication to
what the integration of $M$ should be.  On the other hand, this
sequence shows that the symplectic groupoid of $M$ is an extension
of the monodromy groupoid $\pi_1(\F)$ by the structure groupoid
$\S(M)$. This suggests a different strategy to integrate a regular
Poisson manifold: assuming $\S(M)$ to be a Lie groupoid one looks
for an extension of $\pi_1(\F)$ by $\S(M)$. In \cite{AlHe} it is
shown that this strategy works under some additional technical
assumptions, and Theorem \ref{structure} is actually an
improvement over the results of \cite{AlHe}.

We now turn to the proofs of the results we have stated so far.
Recall that the normal bundle $\nu$ (hence also any associated
tensor bundle) has a natural flat $\F$-connection $\nabla:
\Gamma(T\F)\times \Gamma(\nu)\to\Gamma(\nu)$, given by
\[ \nabla_{X} \overline{Y}= \overline{[X, Y]}.\]
In terms of the Lie algebroid $T\F$, $\nabla$ is a flat Lie
algebroid connection (see \cite{Fer1}) giving a Lie algebroid
\textbf{representation} of $T\F$ on the vector bundle $\nu$ (and
also, on any associated tensor bundle). Hence, the foliated forms
with coefficients in $\nu$, denoted $\Omega^\bullet(\F;\nu)=
\Gamma(\wedge^\bullet T^*\F\otimes \nu)$, carry a foliated de Rham
operator
\begin{align*}
d_{\F}\omega(X_1,\dots,X_{p+1})=&
\sum_{i}(-1)^{i}\nabla_{X_i}(\omega(X_1,\dots,\hat{X_i},\dots,X_{p+1}))\\
+\sum_{i<j}(-1)&^{i+j-1}\omega([X_i, X_j],
X_1,\dots,\hat{X_i},\dots, \hat{X_j},\dots,X_{p+1})),
\end{align*}
for which the associated cohomology, called \textbf{foliated
cohomology} with coefficients in $\nu$, is denoted $H^*(\F;\nu)$.
Similarly one can talk about cohomology with coefficients in any
tensorial bundle associated to $\nu$. The special case of trivial
coefficients $\mathbb{R}$, will be denoted $H^*(\F)$.

For instance, the Poisson tensor in $M$ determines a foliated
2-form $\omega\in\Omega^2(\F)$, which is just another way of
looking at the symplectic forms on the leaves. Therefore, we have
a foliated cohomology class in the second foliated cohomology
group:
\[ [\omega]\in H^2(\F).\]
On the other hand, we saw above in ``Description 2'' of the
monodromy groups, that a splitting $\sigma$ determines a foliated
2-form $\Omega_\sigma\in\Omega^2(\F;\nu^*)$, with coefficients in
the co-normal bundle, and that the corresponding foliated
cohomology class
\[ [\Omega]\in H^2(\F; \nu^*) \]
does not depend on the choice of splitting.

These two classes are related in a very simple way. In fact, there
is a map
\[ d_{\nu}: H^2(\F)\rmap H^2(\F; \nu^*),\]
which can be described as follows. We start with a class
$[\theta]\in H^2(\F)$, represented by a foliated 2-form $\theta$.
As with any foliated form, we have $\theta=
\widetilde{\theta}|_{T\F}$ for some 2-form
$\widetilde{\theta}\in\Omega^2(M)$. Since
$d\widetilde{\theta}|_{\F}= 0$, it follows that the map
$\Gamma(\wedge^2\F)\to\Gamma(\nu^*)$ defined by
\[ (X,Y)\mapsto d\widetilde{\theta}(X, Y, -), \]
gives a closed foliated 2-form with coefficients in $\nu^*$. It is
easily seen that its cohomology class does not depend on the
choice of $\widetilde{\theta}$, and this defines $d_{\nu}$. Now
the formula $\omega(X_{f},X)= X(f)$, immediately implies that
\[ d_{\nu}([\omega])= [\Omega].\]

We emphasize that the operator $d_{\nu}$ is well known in
foliation theory and is part, together with $d_{\F}$, of the
spectral sequence of a foliation. The advantage of this point of
view comes from the fact that the construction of $d_{\nu}$ is
functorial with respect to foliated maps (i.e.~maps between
foliated spaces which map leaves into leaves).

{From} this perspective, a deformation $\gamma_{t}$ of 2-spheres is
a foliated map $S^2\times I\to M$, where in $S^2\times I$ we
consider the foliation $\F_{0}$ whose leaves are the spheres
$S^2\times\{t\}$. {From} $H^2(S^2)\simeq \mathbb{R}$, we get
\begin{align*}
H^{2}(\F_{0})&\simeq C^{\infty}(I)=\Omega^0(I)\\
H^{2}(\F_{0}; \nu)&\simeq C^{\infty}(I) dt=\Omega^1(I),
\end{align*}
where the isomorphisms are obtained by integrating over $S^2$.
Hence, $d_{\nu}$ for $\F_{0}$ becomes the de Rham differential
$d:\Omega^0(I)\to\Omega^1(I)$. Now, the functoriality of $d_{\nu}$
with respect to $\gamma_{t}$ when applied to $\omega$ gives
\[ \frac{d}{dt} \int_{S^2}\gamma_{t}^{*}\omega=
\langle \int_{\gamma_{t}}\Omega,
\frac{d}{dt}\gamma_{t}(p_{N})\rangle.\] This proves the last part
of the proposition (and also the properties of the variation of
the symplectic area, stated at the beginning of this section).

Next, we recall (see \cite{AlHe}) that the second homotopy groups of
the leaves fit into a smooth \'etale groupoid (actually, a bundle of
Lie groups)
\[ \pi_{2}(\F)= \bigcup_{x\in M} \pi_{2}(L_{x}, x) .\]
Also, according to ``Description 3'' in Section \ref{monodromy},
we should view $\partial$ as part of a long exact sequence (but
now of groupoids rather than groups):
\[
\xymatrix{\cdots\ar[r]&\pi_{2}(\F)\ar[r]^{\partial}& \nu^*(\F)
\ar[r]& \Sigma(M)\ar[r]^{\Phi}& \pi_1(\F)}.
\]
This clearly implies the exactness of the sequence of Theorem
\ref{structure}. Now, the smoothness of $\partial$, together with
the fact that $\pi_2(\F)$ is \'etale, imply the following property
of the monodromy groups: for any $a\in \NN_{x}$, there is a smooth
local section $\al$ of $\nu^*(\F)$ defined on some open $U$
containing $x$ such that $\al(x)=a$, and $\al(y)\in \NN_{y}$ for
all $y\in U$. This immediately implies that $r_\NN$ is lower
semi-continuous.

Let us turn now to the equivalence of the various statements in
Theorem \ref{structure}.  We know from the general Lie algebroid
case that (i) implies (ii). The fact that (ii) implies (iii)
follows from the exact sequence above. To prove that (iii) implies
(iv), we just have to observe that in the short exact sequence
\[
\xymatrix{1\ar[r]&\AAA^{'}(M)\ar[r]& \nu^*(\F)\ar[r]& \S(M)\ar[r]&
1},
\]
the projection $\nu^*(\F)\to \S(M)$ will be a submersion if $\S(M)$ is
smooth. We are left with showing that (iv) implies (i), and for that
we will check the integrability conditions of Theorem
\ref{strong-form}. First of all, (iv) implies that $\NN_{x}$ are closed
subgroups of $\nu^{*}_{x}(\F)$. But $\NN_{x}= \partial(\pi_{2}(L_{x},
x))$ are at most countable (because second homotopy groups of
manifolds are so), hence $\NN_{x}$ must be discrete. This shows that
$\AAA^{'}(M)$ must be \'etale, which, in turn, implies that
$\liminf_{y\to x}r_{\NN}(y)\neq 0$. If not, we can find a sequence
$a_{n}\in \NN_{x_{n}}$ of non-zero elements, with $x_{n}\to x$,
$a_{n}\to 0_{x}$. This is impossible since we can find a neighborhood
of $0_{x}$ which only contains zero elements. This concludes our
proof.

\begin{remark}
The last part of this proof can be restated as an interesting
property of any integrable, regular, Poisson manifold: every
transverse deformation $\var(\gamma_t)$ of a sphere $\gamma_0$
such that $\partial([\gamma_0])=0$ must be a trivial deformation,
i.~e., $\partial([\gamma_t])=0$ for all small enough $t$. By the
results of this section, we can rewrite this as:
\[ A'_\omega(\gamma_0)=0\quad \Longrightarrow\quad
A'_\omega(\gamma_t)=0, \ \forall t.\] This property plays an
important role in Alcade Cuesta and Hector approach to
integrability of regular Poisson manifolds, who state it as
``every symplectic vanishing deformation is trivial'' (see
\cite{AlHe}).
\end{remark}

%%%%%%%%%%%%%%%%%%%%%%%%%%%%%%%%%%%%%%%%%%%%%%%%%%%
%%%%%%%%%%%%%%%%%%%%%%%%%%%%%%%%%%%%%%%%%%%%%%%%%%%
\section{Examples}%
\label{Integrability:general}                     %
%%%%%%%%%%%%%%%%%%%%%%%%%%%%%%%%%%%%%%%%%%%%%%%%%%%
%%%%%%%%%%%%%%%%%%%%%%%%%%%%%%%%%%%%%%%%%%%%%%%%%%%

In this section we present several examples of integrable and
non-integrable Poisson manifolds. First, we give some immediate
consequences of Theorem \ref{strong-form} and of the main
properties of the monodromy groups.

\begin{corollary}
If all the monodromy groups $\NN_{x}$ of the Poisson manifold
vanish, then $M$ is integrable. This happens for instance if (i)
for any $x\in M$, the symplectic leaf $L$ through $x$ has finite
second homotopy groups; or if (ii) the isotropy Lie algebra
$\nu^{*}_{x}(L)$ is semi-simple.
\end{corollary}

Not every Poisson manifold is locally integrable (see the example in
Section \ref{example:basic} below). However, from the local form of regular foliations
we deduce:\footnote{As pointed out by Alan Weinstein in \cite{}, 
this corollary should be attributed to Lie, who found the normal
coordinates near a regular point and used them to prove a version 
of ``Lie's third theorem''.}

\begin{corollary}[Lie]
If $x\in M$ is regular, then $M$ is locally integrable around $x$,
i.e., there exists a neighborhood of $x$ in $M$ which is
integrable.
\end{corollary}

After these simple criteria of integrability, we now turn to the
examples. \vskip 15 pt

\subsection{Poisson manifolds of dimension 2}
\label{2-dim} 
The lowest dimension one can have non-trivial Poisson manifolds is
2. However, it follows immediately from Theorem \ref{strong-form} that
in dimension 2 all Poisson manifolds are integrable.

\begin{corollary}
\label{cor:2-dim} Any 2-dimensional Poisson manifold is
integrable.
\end{corollary}

In particular, we recover the main positive integrability result of
\cite{CaFe} that states that any Poisson structure on $\mathbb{R}^2$
is integrable.

Corollary \ref{cor:2-dim} can be partially generalized to higher dimensions
in the following sense: any $2n$-dimensional Poisson manifold whose
Poisson tensor has rank $2n$ on a dense, open set, is
integrable. This, in turn, is a consequence of a general result due to
Debord \cite{Deb} (see also \cite{CrFe}, Corollary 5.9) that states
that a Lie algebroid with almost injective anchor is integrable. The
proof is more involved.

\subsection{A non-integrable Poisson manifold}
\label{example:basic} 
Already in dimension 3 there are examples of non-integrable Poisson manifolds.
We consider $M=\mathbb{R}^3$ with the Poisson bracket:
\[ \{f,g\}= \det
\begin{pmatrix}
x & y & z \\
\frac{\partial f}{\partial x}&\frac{\partial f}{\partial y}&
\frac{\partial f}{\partial z}\\
\frac{\partial g}{\partial x}&\frac{\partial g}{\partial y}&
\frac{\partial g}{\partial z}
\end{pmatrix}.
\]
Identifying $\mathbb{R}^3$ with $\mathfrak{su}(2)^*$, this is just
the Kirillov-Poisson bracket (see also the next example).  We
choose any smooth function $a=a(R)$ on $M$, which depends only on
the radius $R$, and which is strictly positive for $R>0$. We
multiply the previous brackets by $a$, and we denote by $M_{a}$
the resulting Poisson manifold. The bracket on $\Omega^1(M_a)$ is
computed using the Leibniz identity and we get
\begin{align*}
[dx^2,dx^3]&=a dx^1+ b x^1 R \bar{n},\\
[dx^3,dx^1]&=a dx^2+ b x^2 R \bar{n},\\
[dx^1,dx^2]&=a dx^3+ b x^3 R \bar{n},
\end{align*}
where $\bar{n}= \frac{1}{R}\sum_i x^idx^i$ and $b(R)= a'(R)/R$.
The bundle map $\#: T^*M_{a}\to TM_{a}$ is just
\[ \#(dx^i)=a \bar{v}^{i}, \qquad i=1,2,3 \]
where $\bar{v}^i$ is the infinitesimal generator of a rotation about the
$i$-axis:
\[
\bar{v}^1=x^3\frac{\partial}{\partial x^2}-x^2\frac{\partial}{\partial
x^3}, \quad \bar{v}^2=x^1\frac{\partial}{\partial
x^3}-x^3\frac{\partial}{\partial x^1}, \quad
\bar{v}^3=x^2\frac{\partial}{\partial x^1}-x^1\frac{\partial}{\partial
x^2}.
\]
The leaves of the symplectic foliation of $M_{a}$ are the spheres
$S_{R}^{2}\subset\mathbb{R}^3$ centered at the origin, and the
origin is the only singular point. To compute the function
$r_{\NN}$, using the obvious metric on $T^*M_{a}$, we restrict to
a leaf $S_{R}^{2}$ with $R>0$, and we use ``Description 2'' of
Section \ref{monodromy}. As splitting of $\#$ we choose the map
defined by
\[ \sigma(\bar{v}^i)= \frac{1}{a}(dx^i-\frac{x^i}{R}\bar{n}),\]
with curvature the center-valued 2-form
\[ \Omega_{\sigma}= \frac{Ra'- a}{a^2R^3}\omega \bar{n},\]
where $\omega= x^1 dx^2\wedge dx^3+ x^2dx^3\wedge dx^1+
x^3dx^1\wedge dx^2$. Since $\int_{S^2_R}\omega=4\pi R^3$ it
follows that
\[
\NN_{(x, y, z)}\simeq 4\pi \frac{Ra'- a}{a^2}\Zz\bar{n}\subset
\Rr\bar{n}.
\]
The canonical generator of $\pi_{2}(S^{2}_{R})$ defines \emph{the}
symplectic area of $S^2_R$, which is easily computable:
\[ A_{a}(R)= 4\pi \frac{R}{a(R)}.\]
We recover in this way the relationship between the monodromy and
the variation of the symplectic area (Proposition \ref{var-area}).
Also,
\[
r_\NN(x,y,z)= \left\{
\begin{array}{ll}
+\infty \qquad&\  \text{if}\ R=0\ \text{ or }  A'_{a}(R)= 0,\\
\\
A'_{a}(R) \qquad&\  \text{otherwise,}
\end{array}
\right.
\]
so the monodromy might vary in a non-trivial fashion, even for nearby
regular leaves. Our computation also gives the isotropy groups
\[
\Sigma(M_{a}, (x, y, z))\cong \left\{
\begin{array}{llll}
\mathbb{R}^3 \qquad&\  \text{if}\ R=0,\ a(0)= 0,\\
SU(2) \qquad&\ \text{if}\ R=0,\ a(0)\neq 0,\\
\mathbb{R} \qquad&\  \text{if}\ R\neq 0,\  A'_{a}(R)= 0,\\
S^1 \qquad&\  \text{if}\ R\neq 0,\ A'_{a}(R)\neq 0
\end{array}
\right.
\]
We leave it to the reader the task of making the description of
the groupoid $\Sigma(M_a)$ more explicit (see, also, \cite{BCWZ},
where this example is further discussed).

\subsection{Heisenberg-Poisson manifolds}
\label{Heisenberg} 
The integrability of the Heisenberg-Poisson manifolds was discussed in
\cite{Wein2}. We recall that the Heisenberg-Poisson manifold $M(S)$
associated to a symplectic manifold $S$, is the manifold $S\times
\mathbb{R}$ with the Poisson structure given by $\{f, g\}= t\{f_{t},
g_{t}\}_S$, where $t$ stands for the real parameter, and $f_{t}$
denotes the function on $S$ obtained from $f$ by fixing the value of
$t$. We can now easily recover the main result of \cite{Wein2}:

\begin{corollary} 
\label{cor:heisenberg}
For a symplectic manifold $S$, the following are equivalent:
\begin{enumerate}
\item[(i)] The Poisson-Heisenberg manifold $M(S)$ is integrable;
\item[(ii)] $\widetilde{S}$ is pre-quantizable.
\end{enumerate}
\end{corollary}

We recall that condition (ii) is usually stated as follows: when
we  pull back the symplectic form $\omega$ on $S$ to a 2-form
$\widetilde{\omega}$ on the covering space $\widetilde{S}$, the
group of periods
\[ \set{\int_\gamma
  \widetilde{\omega}:\gamma\in\H_2(\widetilde{S},\Zz)}\subset \Rr\]
is a multiple of $\mathbb{Z}$. Note that this group coincides with the
\textbf{group of spherical periods} of $\omega$
\[ \P(\omega)=\set{\int_\gamma \omega:\gamma\in\pi_2(S)},\]
so that (ii) says that $\P(\omega)\subset\Rr$ is a multiple of
$\Zz$.

\begin{proof}[Proof of Corollary \ref{cor:heisenberg}]
We have to compute the monodromy groups. The singular symplectic
leaves are the points in $S\times \{0\}$ and they clearly have
vanishing monodromy groups. The regular symplectic leaves are the
submanifolds $S\times \{t\}$, where $t\neq 0$, with symplectic
form $\omega/t$. The most straightforward way to compute the
monodromy groups is to invoke Proposition \ref{var-area} of
Section \ref{Integrability:regular} (as an instructive exercise,
the reader may also try to use ``Description 2''of Section
\ref{monodromy}). We immediately get
\[ \NN_{(x, t)}= \frac{1}{t}\P(\omega)\subset \mathbb{R},\]
and the result follows.
\end{proof}

Note that one can be quite explicit in describing $\Sigma(M(S))$.
Straightforward computations show that it consists of two types of arrows:
\begin{enumerate}
\item[(a)] arrows which start and end at $(x, 0)$, which form a group
isomorphic to the additive subgroup of $T_{x, 0}^{*}M(S)$;
\item[(b)] arrows inside the symplectic leaves $S\times \{t\}$, $t\neq
0$, which consist of equivalence classes of pairs $(\gamma, v)$, where
$\gamma$ is a path in $S$ and $v\in \mathbb{R}$. Two such pairs
$(\gamma_i,v_i)$ are equivalent if and only if there is a homotopy
$\gamma(\eps,s)$ (with fixed end points) between the $\gamma_{i}$'s,
such that $v_1-v_0=\frac{1}{t}\int\gamma^*\omega$.
\end{enumerate}
This explains how the general strategy of putting a smooth
structure on $\Sigma(M)$ (see the end of Section
\ref{symplectic-groupoid}) is related with the blowing-up
techniques used in \cite{Wein2}.

\begin{remark}
Condition (v) of Theorem \ref{strong-form} splits into two conditions:
\begin{enumerate}
\item[(a)] Each monodromy group $\NN_x$ is discrete, i.e.,
  $r_\NN(x)>0$,
\item[(b)] The monodromy groups are uniformly discrete, i.e.,
  $\liminf_{y\to x}r_\NN(y)>0$.
\end{enumerate}
If we choose a symplectic manifold which does not satisfy the
pre-quantization condition, the Poisson-Heisenberg manifold of Section
\ref{Heisenberg} gives a non-integrable Poisson manifold that violates
condition (a). On the other hand, the example of Section
\ref{example:basic} produces non-integrable Poisson manifolds in which
condition (a) is satisfied, but condition (b) is not. Hence the
two conditions are independent.
\end{remark}

\subsection{Linear Poisson structures}
\label{su(3)} 
Let $\gg$ be any Lie algebra and $G$ a simply connected Lie group with
Lie algebra $\gg$. If we consider the Poisson manifold $M=\gg^*$ with
the Kirillov-Kostant Poisson bracket, and let $A=T^*\gg^*$ we always
obtain an integrable Lie algebroid: the source-simply connected Lie
groupoid integrating $\gg^*$ is
\[ \Sigma(\gg^*)=G\times \gg^* \]
with source and target maps
\[ \s(g,\xi)=\xi,\qquad \t(g,\xi)=\Ad^*g\cdot\xi,\]
and with multiplication
$(g_1,\xi_1)\cdot(g_2\cdot\xi_2)=(g_1g_2,\xi_2)$, wherever
defined. In spite of the fact that linear Poisson structures are
always integrable, the symplectic geometry of their leaves varies in a
non-trivial fashion, and their monodromy reflects this behavior.

For a specific example take $M=\mathfrak{su}^*(3)$. The symplectic
leaves (i.e., the coadjoint orbits) are isospectral sets, and so
we can understand them by looking at their point of intersection
with the diagonal matrices with imaginary eigenvalues. There are
orbits of dimension 6 (distinct eigenvalues), dimension 4 (two
equal eigenvalues) and the origin (all eigenvalues equal). Let us
take for example the (singular) orbit $L$ through
\[ x=\left(\begin{array}{ccc}
i\lambda &0 &0\\ 0&i\lambda&0\\ 0&0&-2i\lambda
\end{array}\right).\]
Then we find its isotropy subalgebra to be
\[
\gg_x=\set{ \left(
\begin{tabular}{c|c}
$\begin{array}{ccc} & &\\ &X&\\ & &
\end{array}$
& 0\\ \cline{1-2} 0& - tr X
\end{tabular}
\right):\ X\in \mathfrak{u}(2)}
\]
and so we see that the simply connected Lie group integrating the
Lie algebra $\gg_x$ is $G(\gg_x)=\Rr\times SU(2)$. We can also
compute the isotropy groups
\begin{align*}
\Sigma(\gg^*, x)&=\set{g\in SU(3): \Ad^*g\cdot x=x}\\ &=\set{
\left(
\begin{tabular}{c|c}
$\begin{array}{ccc} & &\\ &g&\\ & &
\end{array}$
& 0\\ \cline{1-2} 0& $\det g^{-1}$
\end{tabular}
\right):\ g\in U(2)}.
\end{align*}
We conclude that the orbit $L$ is diffeomorphic to
$SU(3)/U(2)=CP(2)$. In fact, one can show that it is
symplectomorphic to $CP(2)$ with a multiple of
its standard symplectic structure
(see \cite{Fer2}, Example 3.4.5). Also, we see that the long exact
sequence
\[ \dots \to \pi_{2}(CP(2),x)\stackrel{\partial}{\to} G(\gg_x)
\to \Sigma(\gg^*, x)\to \pi_1(CP(2),x),\] reduces to:
\[ \dots \to \Zz\stackrel{\partial}{\to} \Rr\times SU(2)
\stackrel{\rho}{\to} U(2)\to \set{1},\] where
$\rho(\theta,A)=e^{i\theta}A$. We conclude that $\partial n=(\pi
n,(-1)^nI)$, so that $\partial$ takes values in the center
$Z(\Rr\times SU(2))=\Rr\times \set{\pm I}$, and
\[ \widetilde{\mathcal{N}}_x=\text{Im~}\partial=2\Zz\times\set{\pm I},\qquad
\mathcal{N}_x=2\Zz.\]
This provides the example promissed in Description \ref{descr2}:
since the last two groups are distinct, there can be no
splitting with center-valued curvature. Another argument is that
such a splitting would define a flat connection on the co-normal
bundle $\nu^*(L)=\Ker\#|_L$, and since $L=CP(2)$ is 1-connected,
it would follow that the co-normal bundle would be a trivial
bundle. This is not possible. In fact, the total Stiefel-Whitney
class of $CP(2)$ is non-trivial, and so is the Stiefel-Whitney
class of the normal bundle, when we embed $CP(2)$ in any Euclidean
space. Hence, the co-normal bundle cannot be trivial.

%%%%%%%%%%%%%%%%%%%%%%%%%%%%%%%%%%%%%%%%%%%%%%%%%%%
%%%%%%%%%%%%%%%%%%%%%%%%%%%%%%%%%%%%%%%%%%%%%%%%%%%
\section{Symplectic realizations}                 %
\label{realizations}                              %
%%%%%%%%%%%%%%%%%%%%%%%%%%%%%%%%%%%%%%%%%%%%%%%%%%%
%%%%%%%%%%%%%%%%%%%%%%%%%%%%%%%%%%%%%%%%%%%%%%%%%%%

Recall (see \cite{Wein4}) that a \textbf{symplectic realization} of a
Poisson manifold $M$ is a surjective Poisson submersion $\mu:S\to M$,
with connected fibers, from a symplectic manifold $S$ onto
$M$\footnote{The reader should be aware that in the literature one
often does not require $\phi$ to be a submersion or to be surjective,
while and the ``symplectic realizations'' defined above are also called
``full symplectic realizations''.}.
The symplectic manifold $S$ comes equipped with a pair of foliations,
in duality with respect to the symplectic structure: the one given by
the fibers of $\mu$, which will be denoted by $\F({\mu})$, and its
symplectic orthogonal $\F(\mu)^{\perp}$. In terms of their tangent
bundles (vectors tangent to the leaves) $\F(\mu)$ is the kernel of the
differential of $\mu$, while $\F(\mu)^{\perp}$ is spanned by the
Hamiltonian vectors $X_{\mu^{*}f}$ ($f\in C^{\infty}(M)$). As
explained in Remark \ref{rem:Hausdorff}, $S$ may be non-Hausdorff, but
we do require the leaves of $\F(\mu)$ and of $\F(\mu)^{\perp}$
to be Hausdorff. In particular, it makes sense to talk about the
completeness of the vector fields $X_{\mu^{*}f}$. The symplectic
realization is called \textbf{complete} if, for any complete Hamiltonian
vector field $X_{f}$ on $M$, the vector field $X_{\mu^{*}(f)}$ is
complete.

The existence of symplectic realizations is guaranteed by the
following basic result:

\begin{theorem}[Karasev \cite{Kar}, Weinstein \cite{CDW}]
Any Poisson manifold $M$ admits a Hausdorff symplectic realization
$\phi:S\to M$.
\end{theorem}

\begin{proof}
Let $M$ be a Poisson manifold. The Lie algebroid $T^*M$ integrates
to a local Lie groupoid $\Sigma_{\text{loc}}(M)$. The proof of
this result in \cite{CaFe} (cf.~Corollary 5.1) shows that one can
exhibit this local Lie groupoid as the quotient
\[ \Sigma_{\text{loc}}(M)=O/\sim\]
where $O$ is a open set in the space of cotangent paths $P(T^*M)$.
Hence, a local version of the construction of Section
\ref{symplectic-groupoid}, shows that one can perform a
Marsden-Weinstein reduction on a open set of $\widetilde{P}(T^*M)$ to
obtain the structure of a local symplectic groupoid on
$\Sigma_{\text{loc}}(M)$. The source map (also the target map)
$\s:\Sigma_{\text{loc}}(M)\to M$ gives a symplectic realization.
\end{proof}

The question of existence of \emph{complete} symplectic realizations
is the main topic of this section, and our main result is the following:

\begin{theorem}
\label{thm:symplectic:realizations} 
A Poisson manifold admits a complete symplectic realization if and
only if it is integrable. 
\end{theorem}

The fact that integrability is somehow related to (and implies) the
existence of complete symplectic realizations is well-known:
if $M$ is integrable, the target map $\s:\Sigma(M)\to M$
gives the desired realization. Completeness follows from the remark
that a path $\gamma:I\to M$ is an integral curve of $X_f$ if and only if it is of
the form $\gamma(t)=\s(\widetilde{\gamma}(t))$ where
$\widetilde{\gamma}:I\to S$ is an integral curve of $X_{f\circ\s}$. Since this 
vector field is left invariant, its flow is left invariant, and the usual 
argument for Lie groups works extends to the groupoid case to show that 
$X_{f\circ\s}$ is complete, provided $X_f$ is complete.

The converse implication, which may look surprising at first sight,
appears to be more difficult since there is no obvious way of
constructing an integrating groupoid out of an arbitrary symplectic
realization $S$ (note that $S$ may be quite different from
$\Sigma(M)$). However, we can now take advantage of the fact that we
always have the groupoid $\Sigma(M)$ around, and the only issue is to
prove its smoothness. This is proved by pulling everything back
to $S$, which one can think of as a desingularization of $M$, where the
problem greatly simplifies.

\begin{remark} We shall see that a complete symplectic realization of
  $M$ can be thought of as faithful representation of the Lie
  algebroid $T^*M$. We can use this representation to integrate
  $T^*M$, in a similar way as one does to integrate a Lie algebra to a
  Lie group, with the help of a faithful representation. Unlike the
  Lie algebra case, when such a representation always exists
 (Ado's
  theorem), in our case we have to assume the existence of such a
  representation, i.e., of the complete symplectic realization $\phi:S\to M$.
\end{remark}

\begin{remark} Note that, in light of Theorem \ref{strong-form} (v),
  the result above can be stated with no reference to integrability or
  groupoids (and this suggests that a different approach might be
  possible).
\end{remark}

\begin{proof}[Proof of Theorem \ref{thm:symplectic:realizations}] 
We let $\mu: S\to M$ be a complete symplectic realization of $M$,
and we first assume that $S$ is Hausdorff. We split the proof into
several steps.
\vskip 5 pt

\emph{Step 1. The Lie algebroid $T^*M$ acts on $S$.} 

The map $df\mapsto X_{\mu^*f}$ defines a bundle map $\rho:
\mu^*T^*M\to TS$, which can also be described as the composition
of the dual of the differential of $\pi$ with the anchor map
$T^*S\to TS$ induced by the symplectic form on $S$. It induces a Lie
algebra homomorphism $\Omega^1(M)\to \X(S)$, and
\begin{equation}
\label{eq-mu} (d\mu)_{s}\circ \rho_{y}= \#_{\mu(y)}, \ \
\mbox{for\ all} \ y\in S .
\end{equation}
One should view $\rho$ as an infinitesimal action of $T^*M$ on $S$
(see also below). Alternatively, we can also think of $\rho$ as the
horizontal lift of a flat, non-linear, contravariant connection (see
\cite{Fer2}). 
\vskip 5 pt

\emph{Step 2. The action of $T^*M$ on $S$ integrates to an action of
  $\Sigma(M)$.}
 
Notice that it is for this reason that we need $\mu$ to be complete (compare
with the integrability of infinitesimal actions of Lie algebras on
manifolds). Given a cotangent path $a:I\to T^*M$, with
base path $\gamma$ starting at $x_{0}\in M$, and given $y\in S$, the
horizontal curve over $a$ starting at $y_{0}$ is the
solution $u: I\to S$ of the initial value problem:
\[ \frac{d}{dt} u(t)= \rho_{u(t)}(a(t)), \ \ u(0)= y_{0} .\]
Let us check that this equation has a unique solution, defined on the
entire unit interval. We choose a time-dependent, compactly supported,
one-form $\alpha$ with the property that $\alpha(t, \gamma(t))= a(t)$,
and we consider the induced time-dependent vector field on $S$, $X(t,
y)= \rho_{y}(\alpha(t, \mu(y)))$. A solution of the equation above is
an integral curve of $X$ with initial condition $y_{0}$, hence
uniqueness. Conversely, if $u$ is an integral curve of $X$,
(\ref{eq-mu}) implies that $\mu\circ u$ is an integral curve of
$\#\alpha$ starting at $x_0$, hence it must be $\gamma$. It
follows that $\frac{d}{dt} u(t)= \rho_{u(t)}(\alpha(t, \gamma(t)))=
\rho_{u(t)}(a(t))$. Finally, $u$ is defined on the entire $I$ because
the completeness of $\mu$ implies that $\rho(\alpha)$ is complete
whenever $\alpha$ is compactly supported.

Next we need to understand how the horizontal lifts depend on the
cotangent path:

\begin{lemma}
\label{lem:homotopy}
Let $a_0$ and $a_1$ be cotangent paths which are cotangent
homotopic. Then their horizontal lifts $u_0$ and $u_1$  are homotopic
paths relative to the end-points.
\end{lemma}

\begin{proof}[Proof of Lemma \ref{lem:homotopy}]
Let us fix the initial point, and consider a family
$a_{\epsilon}$ of cotangent paths, so that $\gamma_{\epsilon}=
\mu\circ a_{\epsilon}$ all start at $x_0$ and end at the same
point. We obtain a corresponding family $u_{\epsilon}$ of paths in
$S$ as above, all starting at $y_{0}$, and all staying in the leaf
of $\F(\mu)^{\perp}$ through $y_{0}$. Moreover,
\begin{equation}
\label{eq-var} \rho (\var(a_{\epsilon}))= \frac{d}{d\epsilon}
u_{\epsilon}(1) ,
\end{equation}
for all $\epsilon$. Of course, this is related to the very
definition of $\var(a_{\epsilon})$: by formula (2) in
Proposition 1.3 of \cite{CrFe}, we have
\begin{equation}
\label{variation}
\var(a_{\epsilon})= \int_{0}^1\phi_{\alpha_{\epsilon}}^{t, s}
\frac{d \alpha_{\epsilon}}{d\epsilon} (s, \gamma_{\epsilon}(s)) ds,
\end{equation}
where $\alpha_{\epsilon}$ are time dependent 1-forms such that
$\alpha_{\epsilon}(t, \gamma_{\epsilon}(t))= a_{\epsilon}(t)$.  Here
$\phi_{\alpha_{\epsilon}}^{t, s}$ denotes the flows of
$\alpha_{\epsilon}$ (see the footnote in
Section \ref{symplectic-groupoid}), which are
covered by the flows of the vector field
$X_\epsilon=\rho(\alpha_{\epsilon})$, and so we find
\[ \rho (\var(a_{\epsilon}))=\int_{0}^1\phi_{X_{\epsilon}}^{t, s}
\frac{d X_{\epsilon}}{d\epsilon} (s,u_{\epsilon}(s)) ds.\]
The left side of this equation is just the variation of parameters
formula for flows of time-dependent vector fields, so (\ref{eq-var})
holds. 

In particular, if $a_{\epsilon}$ is a cotangent homotopy we have
$\var(a_{\epsilon})=0$ and so the end-points $u_{\epsilon}(1)$ are
fixed.
\end{proof}

Hence, for all $g= [a]\in \Sigma(M)$ with $\s(g)=x_{0}$ and
$\t(g)=x_{1}$, given $y_{0}\in S$ in the fiber above $x_{0}$, there is
a well-defined element $u(1)\in S$ in the fiber above $x_{1}$, and we
set $g\cdot y_{0}\equiv u(1)$. In other words we have an induced
action of $\Sigma(M)$ on $S$ (the operation $g\cdot y$ does satisfy
the usual axioms of an action).
\vskip 5 pt

Let us denote by $\mu^{*}\Sigma(M)$ the space of pairs $(g, y)$, with
$\mu(y)= \s(g)$. Then $\mu^{*}\Sigma(M)$ is a groupoid over $S$,
called the \textbf{action groupoid} associated with the action of
$\Sigma(M)$ on $S$: the source and target maps are given by $\s(g,y)=y$ and
$\t(g,y)=g\cdot y$, while the multiplication is given by
\[ (g,y)\cdot (h,z)\equiv (g\cdot h,z),\quad \text{whenever }y=h\cdot z.\]
\emph{Step 3. The action groupoid of $\Sigma(M)$ on $S$ is isomorphic
  to the homotopy groupoid of $\F(\mu)^{\perp}$.} Note that equation (\ref{eq-var}) above actually shows that
$a_{\epsilon}$ is a cotangent homotopy if and only if
$u_{\epsilon}$ is a homotopy with fixed end points. Recall that,
similar to $\Sigma(M)$, one has a homotopy groupoid $G(\F)$
associated to any regular foliation $\F$: it consists of homotopy
classes of paths with fixed end-points, staying in a single leaf.
Now, since $\rho$ is actually an isomorphism from $\mu^*T^*M$ to
the symplectic orthogonal foliation $\F(\mu)^{\perp}$,
(\ref{eq-var}) shows that
\begin{equation}
\label{isom-up} \mu^*\Sigma(M)\cong G(\F(\mu)^{\perp})
\end{equation}
an isomorphism of topological groupoids. 
\vskip 5 pt

\emph{Step 4. $\Sigma(M)$ is a Lie groupoid.} 

Since the homotopy groupoid $G(\F)$ of any regular foliation is always
smooth, the action groupoid $\mu^*\Sigma(M)$ is also smooth and we
want to conclude from this that $\Sigma(M)$ must be smooth too. To
prove it, we will use the exponential map $\exp_{\nabla}: T^*M\to
\Sigma(M)$ with respect to a connection $\nabla$ on $T^*M$, and verify
condition (iv) of Theorem \ref{strong-form}.

Consider the induced map $\mu^{*}T^*M\to \mu^{*}\Sigma(M)$, denoted
$\mu^*\exp_{\nabla}$, and also its composition $F: \mu^{*}T^*M\to
G(\F(\mu)^{\perp})$ with the homeomorphism (\ref{isom-up}). Then $F$
associates to a pair $(v, y)$ ($y\in S$, $v\in T_{\mu(y)}^{*}M$) the
homotopy class of the path $t\mapsto \exp_{\nabla}(tv)y$. This map is a
local diffeomorphism at the zero section. Indeed, after the
identification $\mu^{*}T^*M\cong \F(\mu)^{\perp}$, it is just the
exponential map of $\F(\mu)^{\perp}$ with respect to the pull-back
connection $\mu^*\nabla$. Alternatively, its differential at a point
$(0, y)\in \mu^{*}T^*M$ is an isomorphism since it fits into a
commutative diagram with exact rows:
\[
\xymatrix{ 0 \ar[r] & T^*M \ar[r]\ar[d]_{\rho_{y}} & T_{(0, y)}
\mu^*T^*M \ar[r]\ar[d]_{dF} & T_{y}S\ar[d]_{id}\ar[r] & 0\\
0 \ar[r] & \F(\mu)^{\perp} \ar[r] & T_{1_{y}}G(\F(\mu)^{\perp})
\ar[r]_-{d\s} & T_{y}S \ar[r] & 0 }
\]
It then follows that $\mu^*\exp_{\nabla}$ is locally injective
around the zero section, which immediately implies the similar
property for $\exp_{\nabla}$. Hence, by (iv) of Theorem
\ref{strong-form}, it follows that $M$ must be integrable.

When $S$ is not necessarily Hausdorff, basically the same proof
applies. All the complete vector fields we have used were actually
tangent to the leaves of $\F(\mu)^{\perp}$, which are Hausdorff. Of
course, we also need to know that the construction of the homotopy
groupoid $G(\F)$ works well (i.e., is a smooth manifold) for any
regular foliation $\F$ (on a possibly non-Hausdorff manifold) with
Hausdorff leaves, but this works exactly as in the Hausdorff case.
\end{proof}

The previous arguments become even more natural when using the
language of algebroids. Recall that an \textbf{action of an algebroid} $A$ over
$M$ on a map $\mu: S\to M$ consists of a bundle map $\rho:
\mu^*A\to TS$ with the property that it induces a Lie algebra
homomorphism $\Gamma(A)\to \X(S)$ and satisfies:
\[ (d\mu)_{s}\circ \rho_{s}= \#_{\mu(s)},\quad\text{ for all }s\in S.\] 
The action is called \textbf{complete} if $\rho(\alpha)$ is a complete
vector field whenever $\alpha\in \Gamma(A)$ has compact support. 

A Lie algebroid action $\rho$ determines a Lie algebroid structure on
the pull-back $\mu^*A$, called the \textbf{action Lie algebroid} of
$\rho$: the anchor is simply $\rho$, while the bracket is uniquely
determined by the Leibniz rule and $[\mu^*\alpha, \mu^*\beta]=
\mu^*[\alpha, \beta ]$. Similarly, if a groupoid $G$ acts on $\mu:S\to
M$, there is an induced groupoid $\mu^*G$ over $S$, called the
\textbf{action Lie groupoid}: its arrows consist of pairs $(g, y)\in
G\times S$ with the property that $\mu (y) = \s(g)$. The source of $(g,y)$ is 
$y$, its target is $g\cdot y$, and the multiplication (product)
$(g,y)\cdot(h,z)$, defined when $y=h\cdot z$, equals
to $(gh,z)$.

The same arguments as in the previous proof show that:
\begin{enumerate}
\item[(a)] A complete action of a Lie algebroid $A$ on $\mu: S\to M$
  determines a (topological) action of the Weinstein groupoid $G(A)$
  on $S$, and $\mu^*G(A)\cong G(\mu^*A)$, as groupoids.  Moreover, for
  any connection $\nabla$ on $A$, the pull-back to $S$ of the
  associated exponential map $\exp_{\nabla}\: A\to G(A)$, identifies
  with the exponential $\exp_{\mu^*\nabla}: \mu^*A\to G(\mu^*A)$ with
  respect to the pull-back connection $\mu^*\nabla$.
\item[(b)] If a Lie algebroid $A$ admits a complete action on $\mu:
  S\to M$ such that the action Lie algebroid $\mu^* A$ is
  integrable, then $A$ is integrable. 
\end{enumerate}
The proof above consisted in observing that any symplectic realization
$\mu:S\to M$ comes equipped with an action $\rho_{\mu}$ of $T^*M$,
which is complete if and only if $\mu$ is complete. Moreover,
$df\mapsto X_{\mu^*f}$ defines a Lie algebroid isomorphism between
$\mu^*T^*M$ and $T\F(\mu)^{\perp}$. Since the latter is integrable,
$T^*M$ must also be integrable.

Let us point out several consequences.

\begin{corollary} 
\label{cor:action}
Any complete symplectic realization $\mu: S\to M$ comes equipped with
a (smooth) locally free action of $\Sigma(M)$ on $S$. The action is
free if and only if the leaves of the symplectic orthogonal
foliation $\F(\mu)^{\perp}$ are simply connected.
\end{corollary}

\begin{proof}
That the action is locally free, but not necessarily free, comes
from the fact that $\mu^*\Sigma(M)$ is isomorphic only to the
homotopy groupoid $G(\F(\mu)^{\perp})$. Freeness corresponds to
the case where $G(\F(\mu)^{\perp})$ is a subset of $S\times S$,
or, equivalently, to the simply connectedness of the
leaves. \end{proof}

In \cite{MiWe}, Mikami and Weinstein have explained that ``symplectic reduction''
can be performed in the general context of actions of symplectic groupoids. With our
``infinitesimal point of view'', we find that all is needed is a complete symplectic realization:

In particular, we find that ``symplectic reduction'' can be performed in the
general context of complete symplectic realizations (compare with the similar
result of Mikami and Weinstein \cite{MiWe}, where one starts with an action of the 
symplectic groupoid):

\begin{corollary} 
For any complete symplectic realization $\mu: S\to M$:
\begin{enumerate}
\item[(i)] there is a natural action of the isotropy group $\Sigma(M,x)$ on
  the fiber $\mu^{-1}(x)$, 
\item[(ii)] if $\mu^{-1}(x)/\Sigma(M,x)$ is smooth, then it carries a
  natural symplectic structure, 
\item[(iii)] every cotangent path $a$, starting at $x$ and
  ending at $y$, induces a bijection $\mu^{-1}(x)/\Sigma(M,
  x)\to\mu^{-1}(y)/\Sigma(M, y)$. It only depends on the homotopy
  class of $a$, and, in the smooth case, it is a symplectic
  diffeomorphism. 
\end{enumerate}
\end{corollary}

\begin{proof}
The reduced symplectic structures comes from the fact that, since the
action of $\Sigma(M)$ on $S$ is locally free, the tangent space of
$\mu^{-1}(x)/\Sigma(M, x)$ at some $y\in \mu^{-1}(x)$ is given by
\[ T_y\F(\mu)^{\perp}_{y}/T_y\F(\mu)\cap T_y\F(\mu)^{\perp}.\] 
Then linear symplectic reduction, shows that we have a non-degenerated
bilinear form induced by the symplectic form on $S$. On the other
hand, the action of the cotangent paths on the symplectic quotients
$\mu^{-1}(x)/\Sigma(M,x)$ is just the one induced from the action of
$\Sigma(M)$ on $S$.
\end{proof}

\begin{corollary} 
\label{cor:dual:pair}
Let $\mu: S\to M$ be a complete symplectic realization, and assume
that the orbit space $S/\Sigma(M)$ is smooth. Then it carries a
natural Poisson structure, whose symplectic leaves can be identified
with the symplectic manifolds $\mu^{-1}(x)/\Sigma(M,x)$.
\end{corollary}

These results show that the map $\mu:S\to M$ can be
can be seen as a \textbf{moment map} for the groupoid action of $\Sigma(M)$ on
$S$. In fact, it is easy to see that the graph of the action
$\set{(g,y,g\cdot y):\s(g)=\mu(y)}$ is a Lagrangian submanifold of
$\Sigma(M)\times S\times\bar{S}$, so the groupoid action of $\Sigma(M)$ on
$S$ is symplectic in the sense of Mikami and Weinstein \cite{MiWe}. In
fact, the corollaries above also follow from this observation hence
the
general results of \cite{MiWe}, since we know now that $\Sigma(M)$ is
smooth (see also the recent preprint \cite{Wein5}).  

Notice that a complete symplectic realization gives rise to a two leg
diagram
\[
\xymatrix{
&S\ar[dl]_{\mu}\ar[dr]^{\pi}\\
M& & S/\Sigma(M)}
\]
where the quotient $S/\Sigma(M)$ and the map $\pi$ are Poisson under
some favorable circumstances (e.g. when
  $S/\Sigma(M)$ is smooth and $S$ is a principal $\Sigma(M)$-bundle over $S/\Sigma(M)$). Moreover, if the fibers of $\mu$ are
simply connected, the symplectic groupoid of
this new Poisson manifold is just the gauge groupoid
$S\times_{\Sigma(M)}S$. Note also that if one repeats this
construction on
$S/\Sigma(M)$ and the symplectic realization given by
$\pi$, then one gets back $M$.

\begin{example}
\label{red-Poisson}
A well-known instance of this is provided by a complete symplectic
realization $J:S\to\gg^*$ of the dual of some Lie algebra. In this
case we have $\Sigma(\gg^*)=\gg^*\times G$, where $G$ is a simply
connected Lie group with Lie algebra $\gg$, and
$\Sigma(\gg^*,\nu)=G_\nu$ is the isotropy group of $\nu\in \gg^*$ for
the coadjoint action (see Section \ref{example:basic}). Then we obtain
a free action of $G$ on $S$ which is Hamiltonian with momentum map
$J$, and the symplectic structures on $J^{-1}(\nu)/\Sigma(\gg^*,\nu)$
are the well-known Marsden-Weinstein symplectic quotients. They form
the leaves of the Poisson manifold $S/\Sigma(\gg^*)=S/G$, provided
this quotient is smooth. We will came back to dual pairs later in our
discussion of Morita equivalence (see Section \ref{Morita}).
\end{example}

Corollaries \ref{cor:action} through \ref{cor:dual:pair} improve and
explain the results of Mikami and Weinstein \cite{MiWe,Wein5}.  We can
summarize this section by saying that we can view a symplectic
realization of a Poisson manifold $M$ as a momentum map for an action
of the groupoid $\Sigma(M)$ on $S$. Being the target of a momentum
map, $M$ can be thought of as the dual of the Lie algebra(oid) of its
Weinstein group(oid). This fits well with Alan Weinstein's
remark (see \cite{CaWe}, page 46) that ``it is tempting to think of
any symplectic manifold $S$
as the dual of the Lie algebra of
$\pi_1(S)$''.

%%%%%%%%%%%%%%%%%%%%%%%%%%%%%%%%%%%%%%%%%%%%%%%%%%%
%%%%%%%%%%%%%%%%%%%%%%%%%%%%%%%%%%%%%%%%%%%%%%%%%%%
\section{Induced Poisson structures}              %
\label{Poisson:category}                          %
%%%%%%%%%%%%%%%%%%%%%%%%%%%%%%%%%%%%%%%%%%%%%%%%%%%
%%%%%%%%%%%%%%%%%%%%%%%%%%%%%%%%%%%%%%%%%%%%%%%%%%%

In this section we discuss submanifolds which have a canonical
induced Poisson structure. 

Let $M$ be a manifold with a smooth foliation $\F$, which may be
singular (for singular foliations, see e.g. Vaisman's book \cite{Vai}). 
By a \textbf{smooth family of symplectic forms} on the
leaves we mean a family of symplectic forms
$\set{\omega_L\in\Omega^2(L):L\in\F}$ such that for every smooth
function $f\in C^\infty(M)$ the Hamiltonian vector field $X_f$ defined
by
\[ i_{X_f}\omega_L=d(f|_L),\quad \forall L\in\F\]
is a smooth vector field in $M$. If $(M,\{\cdot,\cdot\})$ is a Poisson
manifold, then the symplectic foliation with the induced symplectic
forms on the leaves, gives a smooth foliation with a smooth family of
symplectic forms.

Conversely, let $M$ be a manifold with a smooth foliation $\F$,
furnished with a smooth family of symplectic forms on the leaves. Then,
we have a Poisson bracket on $M$ defined by the formula
\[ \{f,g\}\equiv X_f(g),\]
for which the associated symplectic foliation is precisely $\F$. Hence, a
Poisson structure can be defined by its symplectic foliation instead
of the Poisson bracket (see also \cite{Vai}, Theorem 2.14). This
motivates the following definition.

\begin{definition}
\label{def:Dirac}
Let $M$ be a Poisson manifold. A submanifold $N\subset M$ is called a
\textbf{Poisson-Dirac submanifold} if $N$ is a Poisson manifold and 
\begin{enumerate}
\item[(i)] The symplectic foliation of $N$ is $N\cap\F=\set{L\cap
    N:L\in\F}$, 
\item[(ii)] For every leaf $L\in\F$, $L\cap N$ is a symplectic
  submanifold of $L$.
\end{enumerate}
\end{definition}

Let us clarify this definition. First, the intersections $N\cap L$
need not be connected, so in (i) we really mean that: 
\begin{enumerate}
\item[(ia)] $N$ intersects $L$ cleanly, so that $N\cap L$ is a
  submanifold of $N$ and $L$ and $T(N\cap L)=TN\cap TL$,
\item[(ib)] the symplectic leaves of $N$ are the connected components
  of the intersections $N\cap L$. 
\end{enumerate}
Also, condition (ii) means that the symplectic forms on $L\cap N$ are
the pull-backs $i^*\omega_L$, where $i:N\cap L\hookrightarrow L$ is the
inclusion. Then we must have
\begin{equation}
\label{eq:Dirac}
TN\cap \#(TN^0)=\set{0},
\end{equation}
since the left-hand side is the kernel of the pull-back $i^*\omega_L$.

By the remarks above, a submanifold $N$ has at most one Poisson
structure satisfying these properties, and this Poisson structure is
completely determined by the Poisson structure of $M$. Poisson
submanifolds are obvious examples of Poisson-Dirac submanifolds, and we will
see many other examples later. 

% Hence, we see that the
% property of being Poisson-Dirac is a property that a submanifold may (or not)
% enjoy, and does not depend on any aditional structure. In a nutshell,
% Poisson-Dirac submanifolds are the submanifolds of a Poisson manifold that
% have a natural induced Poisson structure.

In order to discuss the problem of integrability of Poisson-Dirac
submanifolds, it will be convenient to recall a few facts about
Dirac structures. At the same time, this will allow us to justify our
usage of the term \textbf{Poisson-Dirac submanifold}. 

\subsection{Poisson-Dirac subspaces}
\label{Linear case}
Recall (see \cite{Cor}) that a \textbf{linear Dirac structure} on a vector
space $V$ is a subspace $L\subset V\oplus V^*$ which is maximally
isotropic with respect to the canonical symmetric pair $\langle\cdot ,
\cdot \rangle$ given by:
\[ \langle(v, \xi), (w, \eta)\rangle= \frac{1}{2}(\xi(w)+ \eta(v)).\]
We remark that linear Dirac structures can always be restricted to
subspaces: if $W\subset V$ is a subspace, then on $W$ one has the
induced Dirac structure
\[ L_{W}=\set{(v,\xi|_{W}):(v,\xi)\in L,v\in W}. \]

Now a Poisson vector space $(V,\Pi)$ (a vector space $V$ with a
bivector $\Pi\in\wedge^2 V$) is the same as a linear Dirac structure
on $V$, with the property that the projection $L\to V^*$ is bijective.
Namely, a bivector $\Pi$ is completely determined by its graph $L^\Pi\subset
V\oplus V^*$, where
\[ L^{\Pi}\equiv \set{(\#\xi,\xi): \xi\in V^*}.\]
Hence, given a subspace $W\subset V$ one has the induced Dirac structure 
\[ L^\Pi_{W}=\set{(\#\xi,\xi|_{W}):\xi\in V^*,\#\xi\in W}. \]
However, in general, the projection $L^\Pi_{W}\to W^*$ will not be
bijective so this Dirac structure is not defined by some bivector in
$W$. Noticing that the kernel of this projection is precisely $W\cap
\#(W^0)$ we introduce the following definition.

\begin{definition}
Let $(V,\Pi)$ be a Poisson vector space. A subspace $W\subset V$ is
called a \textbf{Poisson-Dirac subspace} if
\begin{equation}
\label{eq:linear:Dirac}
W\cap \#(W^0)=\set{0}.
\end{equation}
\end{definition}

Therefore, a Poisson-Dirac subspace $W$ of a Poisson vector space $(V,\Pi)$
has a natural induced bivector $\Pi_W$. Let us set
\[A_W\equiv(\#(W^0))^0=\set{\xi\in V^*:\#\xi\in W}.\]
The bivector $\Pi_{W}$ is given by
\begin{equation}
\label{Dirac:Poisson:tensor}
\Pi_{W}(\xi, \eta)=\Pi(\tilde{\xi}, \tilde{\eta}), 
\end{equation}
where $\tilde{\xi}$ and $\tilde{\eta}$ are extensions to $V$ of
$\xi$ and $\eta$, at least one of which lies in $A_W$. On the other hand,
condition (\ref{eq:linear:Dirac}) can also be written as:
\[ W^0+A_W=T^*M.\]
On the other hand, if (\ref{eq:linear:Dirac}) holds, we see that
\[ W^0\cap A_W=W^0\cap\Ker\#,\]
so we have a short exact exact sequence
\[\xymatrix{0\ar[r]& W^0\cap Ker\# \ar[r]& A_{W}\ar[r]& W^*\ar[r]& 0}.\]
The space $A_{W}$ has the role of relating $V^*$ and $W^*$, and
shows that $\Pi_{W}$ is obtained by linear presymplectic
reduction. Indeed, $A_{W}$ is a subspace of the (presymplectic)
vector space $V^*$, while $W^*$ is a quotient of $A_{W}$ with the
quotient map just the restriction to $W$.

We shall call a \textbf{Dirac projection} of $W$ any projection
$p:V\to W$ such that $p|_{\#(W^0)}= 0$. Notice that $p: V\to W$ is a
Dirac projection if and only if $p^*:W^*\to V^*$ is a splitting of the
short exact sequence above. The role of a Dirac projection $p:V\to W$
is to replace 
the choice of extensions in the description
(\ref{Dirac:Poisson:tensor}) of $\Pi_{W}$:
\[ \Pi_W(\xi,\eta)=\Pi(p^*\xi,p^*\eta).\]
Any projection $p:V\to W$ determines a splitting $V= W\oplus E_{p}$,
where $E_{p}=(1-p)V$, and Dirac projections correspond to
$\Pi$-orthogonal complements of $W$ in $V$ (that is, $\Pi(\xi,\eta)=
0$ for all $\xi\in W^0$, $\eta\in E_{p}^0$). This orthogonality
condition shows that $\Pi$ decomposes as
\[\Pi= \Pi_{W}+ \Pi_{E_p},\]
for a unique $\Pi_{E_p}\in\wedge^2 E_p$.

We define the \textbf{rank} of a Poisson-Dirac subspace $W\subset V$ to be the number
\[ \rank W= \dim (W^0\cap \Ker\#).\]
This number $r=\rank W$ determines the dimensions of the other spaces
involved. For instance, $\dim A_{W}=\dim W+r$ and $\dim V=\dim
(W+\#(W^0))+r$. We find that 
\[ \codim W-\rank W=\rank \Pi -\rank \Pi_N,\]
so that $0\leq \rank W \leq \codim W$.

\begin{example} 
A \textbf{cosymplectic subspace} of a Poisson vector
space $(V,\Pi)$ is a subspace $W\subset V$ such that $V=
W+\#(W^0)$. Hence, a cosymplectic subspace is the same thing as a
Poisson-Dirac subspace of minimal rank (i.e., $\rank W=0$). Since this
is the same as satisfying the transversality condition $W+ Im(\#)= V$,
a cosymplectic subspace is a non-zero Poisson-Dirac subspace which has a
unique Dirac projection.

At the other extreme, we can consider Poisson-Dirac subspaces $W\subset V$ of
maximal rank (i.e., $\rank W=\codim W$). Then we must have
$W^0=\ker\#$ or, equivalently, $\Im\#_{W}=\Im \#$. This just means
that $W$ is a \textbf{Poisson subspace}.
\end{example}

Finally, note that a Poisson vector space $(V,\Pi)$ is completely
determined by the linear symplectic space $(S,\omega)$, where
$S=\Im\#\subset V$ and 
\[ \omega(\#\xi, \#\eta)= \Pi(\xi,\eta).\]
This is the linear counterpart of the remark made at the beginning of
the section, since the leaves of $\Pi$ are just the linear symplectic
affine spaces obtained from $(S,\omega)$ by translation.

Now for any Poisson-Dirac subspace $W\subset V$, we have
\[ \Im \#_{W}= W\cap \Im \#.\]
Hence, the linear symplectic space $(S_W,\omega_W)$ that corresponds
to $W$, is such that $S_W\hookrightarrow S$ and $\omega_W=i^*\omega$,
with $i$ the inclusion. Conversely, the canonical splitting
\[ \Im \#= \Im\#_{W}\oplus \#(W^0),\]
shows that if $W\cap Im(\#)$ is a symplectic subspace of $Im(\#)$
then $W$ is Poisson-Dirac. 

Therefore, if we view a Poisson vector space $(V,\Pi)$ as a linear
symplectic space $(S,\omega)$, the Poisson-Dirac subspaces of $V$ correspond
to the symplectic subspaces of $S$.

%%%%%%%%%%%%%%%%%%%%%%%%%%%%%%%%%%%%%%%%%%%%%%%%%%%
\subsection{Poisson-Dirac submanifolds}           %
\label{Global case}
%%%%%%%%%%%%%%%%%%%%%%%%%%%%%%%%%%%%%%%%%%%%%%%%%%%

We now consider the non-linear case. Let $(M,\Pi)$ be a Poisson
manifold, and consider a submanifold $N\subset M$ which is
``pointwise Poisson-Dirac'', i.e., $T_x N$ is a Poisson-Dirac subspace of $T_x M$
(with respect to $\Pi_x\in \wedge^2 T_x M$) for all $x\in N$: 
\[ T_x N\cap \#(T_x N^0)=\set{0}.\]
{From} the linear case it follows that there is an induced two-tensor on $N$,
denoted $\Pi_{N}$, but there is nothing to ensure us that it is
smooth.  

\begin{example}
Let $M=\Cc^3$ with complex coordinates $(x,y,z)$. We consider the
(regular) foliation of $\Cc^3$ by complex lines defined by
\[ dy=0,\quad dz-ydx=0.\]
The leaves of this foliation are symplectic submanifolds of $\Cc^3$
with the canonical symplectic form. Hence, we have a Poisson structure
on $M=\Cc^3$ with this symplectic foliation, and we denote it by
$\Pi$. Now consider the submanifold $N=\set{(x,y,z):z=0}\subset
M$. A complex line in the foliation intersects $N$ in a point (if
$y\neq 0$) 
or in a complex line (if $y=0$), which are symplectic
submanifolds. This shows, that 
\[ TN\cap\#(TN^0)=\set{0},\]
so $N$ is pointwise Poisson-Dirac, 
 and $\Pi_N$ is not smooth, for its
image is a non-smooth distribution (notice that it is, none less an
integrable distribution!).
\end{example}

If $\Pi_N$ is smooth, the
  next proposition shows 
that $\Pi_{N}$ satisfies
automatically the integrability condition $[\Pi_{N},\Pi_{N}]=0$ (i.e.,
the induced almost Dirac structure on $N$ is integrable):

\begin{proposition}
Let $N$ be a submanifold of the Poisson manifold $M$, such that
\begin{enumerate}
\item[(i)] $N$ is pointwise Poisson-Dirac, i.e., $TN\cap\#(TN^0)=\set{0}$,
\item[(ii)] The induced tensor $\Pi_N$ is smooth.
\end{enumerate}
Then $\Pi_N$ is a Poisson tensor on $N$.
\end{proposition}

\begin{proof}
Fix any $x\in N$. We claim that $[\Pi_N,\Pi_N]_x=0$. By the linear
theory, we can choose a splitting $T_N M=T N\oplus E$ where $E$ is
some vector bundle over $N$, such that
\[\Pi=\Pi_N+\Pi_E,\]
where $(\Pi_E)_x\in \wedge^2 E_x$, and $(E_{x})^{0}$ and $(T_x N)^0$ are
$\Pi_x$-orthogonal complements (\footnote{Notice that for
  $y\not= x$, in general, we will have
  $(\Pi_E)_y=(\Pi-\Pi_N)_y\not\in\wedge^2 E_y$, but that is irrelevant
  for the argument.}). 
Since $[\Pi,\Pi]=0$, it follows that 
\begin{align*} 
[\Pi_N,\Pi_N]_x&=[\Pi_N+\Pi,\Pi_N-\Pi]_x,\\
               &=-[2\Pi_N+\Pi_E,\Pi_E]_x.
\end{align*}
The left-hand side of this expression lies in $\wedge^3 T_x N\subset
\wedge^2T_x M\wedge T_xN$, while the right-hand side lies in
$\wedge^2T_x M\wedge E_x$. Hence they must both be zero.
\end{proof}

This proposition has the following corollary which justifies our usage of the
term \textbf{Poisson-Dirac submanifold} in Definition \ref{def:Dirac}.

\begin{corollary}
\label{cor:Dirac}
A submanifold $N\subset M$ of a Poisson manifold is a Poisson-Dirac submanifold
iff it is pointwise Poisson-Dirac and the induced tensor is smooth.
\end{corollary}

\begin{proof}
If $N\subset M$ is a Poisson-Dirac submanifold then it is pointwise
Poisson-Dirac. Also, the induced tensor coincides with the Poisson tensor on
$N$ so it is smooth.

Conversely, if $N\subset M$ is pointwise Poisson-Dirac and the induced tensor
is smooth, then by the proposition $\Pi_N$ is a Poisson tensor, so $N$
is a Poisson manifold. By the linear theory, the smooth, integrable,
distribution defined by $\Pi_N$ on $N$ (namely $\Im(\#_N)$) is just
$TN\cap\Im(\#)$. Hence, the symplectic foliation of $N$ is $N\cap\F$,
where $\F$ is the symplectic foliation of $M$, and the leaves of $N$ are
symplectic submanifolds of the leaves of $M$. Therefore, $N$ is a
Poisson-Dirac submanifold.
\end{proof}

A \textbf{Dirac projection} of $N$ is a smooth bundle map $p:T_N M\to
TN$ with the property that $p_x:T_x M\to T_x N$ is a linear Dirac
projection for each $x\in N$. In other words, $\{p_{x}\}$ is a family
of Dirac projections depending smoothly on $x$. The following
proposition shows that Poisson-Dirac submanifolds which admit a Dirac
projection are the objects which Vaisman calls in \cite{Vai2}
\emph{quasi-Dirac submanifolds}.

\begin{proposition}
\label{prop:Dirac:projection}
Let $M$ be a Poisson manifold and $N\subset M$ a submanifold. The
following statements are equivalent:
\begin{enumerate}
\item[(i)] $N$ is a Poisson-Dirac submanifold admitting a Dirac projection
  $p:T_N M\to TN$; 
\item[(ii)] There exists a bundle $E$ such that $T_N M=TN\oplus
  E$ and $\#(E^0)\subset TN$.
\end{enumerate} 
\end{proposition}

\begin{proof}
To show that (i) $\Rightarrow$ (ii), we just observe that if $p:T_N
M\to TN$ is a Dirac projection then $E=\Ker p$ is a subbundle such
that $T_N M=TN\oplus E$ and $\#(E^0)\subset TN$.

For the converse, we observe (as in
the linear case) that (ii) gives a decomposition
\[ \Pi=\Pi_N+\Pi_E\]
with $\Pi_N\in \Gamma(\wedge^2 TN)$ and $\Pi_E\in \Gamma(\wedge^2
E)$, where now both $\Pi_N$ and $\Pi_E$ must be smooth bivector
fields. By Corollary \ref{cor:Dirac}, $N$ is a Poisson-Dirac submanifold.
\end{proof}

\begin{remark}
If, around each point of $N$, we have local Dirac projections we can
use a partition of unity to glue them into a global Dirac
projection. In fact, notice that a convex combination of linear Dirac
projections is a linear Dirac projection. Therefore, the existence of
a Dirac projection for a given Poisson-Dirac submanifold is a local issue.
\end{remark}

\begin{corollary}
Let $N$ be a submanifold of Poisson manifold which is pointwise Poisson-Dirac
and for which the Poisson-Dirac subspaces $T_x N\subset T_xM$ have constant
rank. Then $N$ is a Poisson-Dirac submanifold admitting a Dirac projection.
\end{corollary}

\begin{proof}
The assumptions imply the existence of a Dirac projection: if $T_x
N\subset T_xM$ have constant rank then $\#(TN^0)$ has constant rank,
and so we can choose a bundle $F$ such that $T_N M=TN\oplus
\#(TN^0)\oplus F$. Then projection onto $TN$ is a Dirac
projection.
\end{proof}

\begin{example}
Let $N$ be a submanifold of a Poisson manifold $M$. If $M$ is a
cosymplectic submanifold (i.e., each $T_{x}N$ is a cosymplectic
subspace of $T_{x}M$), then $N$ is a Poisson-Dirac submanifold.  Similarly, if
$N$ is a Poisson submanifold (i.e., each $T_{x}N$ is a Poisson
subspace of $T_{x}M$), then $N$ is also a Poisson-Dirac submanifold.
\end{example}

For a general Poisson-Dirac submanifold $N\subset M$, the rank of the linear
Poisson-Dirac subspaces $T_x N\subset T_x M$ will vary. Let us
call this number the \textbf{rank of the Poisson-Dirac submanifold} at
$x$, and denote it by $\rank_x N$. By the linear theory, we have:
\[ \codim N-\rank_x N=\rank \Pi_x-\rank (\Pi_N)_x.\]

Obviously, cosymplectic and Poisson submanifolds are examples of 
Poisson-Dirac submanifolds of constant rank. Here is a simple example of a
Poisson-Dirac submanifold with non-constant rank.

\begin{example}
Let $M=\Rr^3$ with coordinates $(x,y,z)$ and the following Poisson
bracket:
\[ \{x,y\}=f(z),\quad \{x,z\}=\{y,z\}=0,\]
where $f\in C^\infty(\Rr)$ is such that $f(z)=0$ for $z\le 0$ and
$f(z)>0$ for $z>0$. Then one checks easily that 
$N=\set{(0,0,z)\in\Rr^3}$ is a Poisson-Dirac submanifold. This Poisson-Dirac submanifold
has rank two for $z\le 0$, and has rank zero for $z>0$. 
\end{example}

\begin{example}
In the previous example, a Dirac projection still exists. Now,
let $M=\Rr^4$ with coordinates $(x,y,z,w)$ and Poisson
bracket defined by:
\[ \{x,y\}=x^2,\quad \{z,w\}=z.\]
Then one checks easily that $N=\set{(x,y,z,w)\in\Rr^4:y=0,x=z^2}$ is
a Poisson-Dirac submanifold. Let us use coordinates $(z,w)$ for $N$. The
induced Poisson bracket on $N$ is simply
\[ \{z,w\}_N=z.\]
and $N$ is a Poisson-Dirac submanifold of rank 2 at points $(0,w)$ and of rank
0 elsewhere. Hence, for $z\not=0$ there is a unique Dirac
projection. The subbundle $E\subset T_N M$, whose existence is
asserted by Proposition \ref{prop:Dirac:projection}, is given by
\[ E_{(z,w)}=\set{(-z^4 a_1,z^4 a_2,0,-2z^2 a_2):a_1,a_2\in\Rr},
\quad (z\not=0).\]
This bundle can be extended over points with $z=0$, where it will have fiber
\[ E_{(0,w)}=\set{(a_1,0,0,a_2):a_1,a_2\in\Rr}.\]
Note, however, that at these points $TN\cap E\not=0$, so in a
neighborhood of any point $(0,w)$ there are no Dirac projections.
\end{example}

Let us restrict now to constant rank Poisson-Dirac submanifolds. Then we have
the vector bundle:
\[ A_N(M)\equiv\set{\xi\in T^*_N M: \#\xi\in TN},\]
and this has in fact a Lie algebroid structure over $N$, with bracket
induced from the bracket on $T^*M$. On the other hand, we
also have the vector bundle
\[ \gg_N(M)\equiv TN^0\cap \Ker\#, \]
which is in fact a bundle of Lie algebras over $N$. Hence, we see that
one has a short exact sequence of Lie algebroids:
\begin{equation}
\label{seq-Dirac}
\xymatrix{0\ar[r] & \gg_N(M) \ar[r] & A_N(M)\ar[r] &T^*N\ar[r] &0}
\end{equation}

The Lie algebroid $A_{N}(M)$ determines a subgroupoid $\Sigma_{N}(M)$
of the Weinstein groupoid $\Sigma(M)$, which consists of equivalence
classes $[a]\in\Sigma(M)$ that can be represented by a cotangent path
$a$ with $\#a\in TN$. The short exact sequence above then gives a
groupoid homomorphism $\Sigma_N(M)\to\Sigma(N)$ with kernel a bundle
of Lie groups $\G(\gg_N)$. We shall show now that, in the integrable
case, the diagram
\[ 
\newdir{ (}{{}*!/-5pt/@^{(}}
\xymatrix{ \Sigma_N(M)\ar@{ (->}[r]\ar[dr]&\Sigma(M)\\& \Sigma(N)}\]
corresponds in fact to symplectic reduction.

\begin{proposition}
Let $N$ be a constant rank Poisson-Dirac submanifold of a Poisson manifold
$M$. If $M$ is integrable then $\Sigma_N(M)$ is a Lie
subgroupoid of $\Sigma(M)$. Moreover, $\Sigma_N(M)\subset\Sigma(M)$ is a
presymplectic submanifold with rank equal to $2\dim N$ and its 
characteristic foliation has leaves the orbits of the equivalence relation on
$\Sigma_N(M)$ defined by the bundle of Lie groups $\G(\gg_N)$.
\end{proposition}

\begin{proof}
Subalgebroids of an integrable Lie algebroid integrate to Lie
subgroupoids of the corresponding Weinstein groupoids. Hence,
the only thing to check is the statement about the characteristic
foliation. By invariance of the symplectic form $\omega\in\Sigma(M)$,
it is enough to check that the kernel of $i^*\omega$, where
$i:\Sigma_N(M)\hookrightarrow\Sigma(M)$ is the inclusion, coincides
with $(\gg_N)_x$ at points $x\in M$.

Now observe that if $x\in M$, we have an isomorphism of symplectic vector
spaces 
\[ T_{x}\Sigma(M)\simeq T_{x}M \oplus T_{x}^{*}M\]
where on the right the symplectic form is the the one defined by
\[ \omega_x((v_1,\xi_1),(v_2,\xi_2))=\xi_1(v_2)-\xi_2(v_1)+\Pi(\xi_1,\xi_2).\]
(the identity section gives a splitting of the differential of the
source map). Under this isomorphism, we have
\[ T_x\Sigma_N(M)\simeq T_x N\oplus A_N(M)_x, \]
and we check that
\[
\omega_x((v_1,\xi_1),(v_2,\xi_2))=0, \quad
\forall v_2\in T_x N,\xi_2\in A_N(M)_x\quad
\Longleftrightarrow \quad v_1=0,\ \xi_1\in\gg_x,
\]
so the result follows.
\end{proof}

We conclude that if both $N$ and $M$ are integrable, the symplectic
form on $\Sigma_N(M)$ is obtained by symplectic reduction from the
symplectic form on $\Sigma(M)$: we first pullback to $\Sigma_N(M)$, to
obtain a presymplectic form, and then we project to $\Sigma(N)$ along
the characteristic foliation. 

\begin{example}
For a cosymplectic submanifold, we have $\rank N=0$ so $\Sigma_N(M)$
is a symplectic subgroupoid of $\Sigma(M)$ of dimension $2\dim N$ 
and we have $\Sigma(N)=\Sigma_N(M)$ (but not the converse; see below).
On the other hand, for a Poisson submanifold we have $\codim N=\rank
N$ and this happens precisely when $\Sigma_N(M)$ is a coisotropic
submanifold of $\Sigma(M)$.
\end{example}

Notice that $M$ being integrable guarantees that $A_N$ is integrable.
Although the foliation of the Lie algebroid $A_N$ coincides with the
symplectic foliation of $N$, this is not enough to guarantee that $N$
is integrable, since the isotropy Lie algebras of $A_N$ and $T^*N$ are
distinct. Also, $\gg_N(M)$ being a bundle of Lie algebras is always
integrable.

\begin{remark}
In general, a Poisson-Dirac submanifold will not have constant rank. However,
  we can still define $\Sigma_N(M)$ which will be a non-smooth
  subgroupoid of $\Sigma(M)$. We will still have a diagram as above
  relating $\Sigma(M)$, $\Sigma(N)$ and $\Sigma_N(M)$. So, morally,
  Poisson-Dirac submanifolds are the base manifolds of presymplectic
  groupoids (see also \cite{BCWZ}), with Poisson submanifolds 
  corresponding to base manifolds of coisotropic subgroupoids.
\end{remark}

One can define a \textbf{presymplectic groupoid} as a Lie groupoid $\G$
with a closed 2-form $\omega\in \Omega^2(\G)$ such that:
\begin{enumerate}
\item[(i)] $\omega$ is compatible with the product:
$m^*\omega=\pi_1^*\omega+\pi_2^*\omega$,
\item[(ii)] $\omega$ has a characteristic foliation defined by a
normal Lie subgroupoid $\H\subset \G$.
\end{enumerate}

Property (i) is similar to the corresponding property for symplectic
groupoids (see Proposition \ref{prop:compatible:symplectic}). Now, to
explain (ii), recall that a normal subgroupoid $\H\subset\G$ is a wide
subgroupoid such that that: 
\[ \forall h\in \H_x, g\in\G \text{ with }s(g)=x\ \Longrightarrow\ 
ghg^{-1}\in \H.\] 
Notice that this condition only involves the isotropy groups of $\H$.
A normal subgroupoid $\H\subset\G$ defines an equivalence relation in
$\G$ and an equivalence relation in the base $N$:
\begin{itemize}
\item For $g_1,g_2\in\G$, $g_1\sim g_2$ if and only if there exist
  $h,h'\in\H$ such that $hg_1h'=g_2$;
\item For $x,y\in M$, $x\sim y$ if and only if there exists $h\in\H$ with
  $s(h)=x$ and $t(h)=y$; 
\end{itemize}
So condition (ii) means that the leaves of the characteristic
foliation are the $\H$ orbits in $\G$.

Presymplectic groupoids are interesting objects which we believe
deserve more attention. Let us list some of the properties we hope
they will satisfy: Let $\G$ be a presymplectic groupoid over $N$ with
characteristic foliation defined by $\H\subset\G$. Then:
\begin{enumerate}
\item[(a)] The base manifold $N$ has a Dirac structure with
null foliation the foliation defined by $H$;

\item[(b)] Conversely, every Dirac structure has a (Weinstein)
  presymplectic groupoid; 

\item[(c)] If the leaf space $N/H$ is smooth it is a Poisson
manifold. Its symplectic groupoid is $\G/\H$ provided it is smooth;

\item[(d)] The groupoid $\G'$ obtained from $\G$ by factoring out only
  the isotropy of $\H$, has Lie algebroid $L_N\subset TN\oplus T^*N$
  (the Dirac subbundle), provided it is smooth;

\item[(e)] By a groupoid-equivariant coisotropic embedding theorem, $\G$
should embed into a symplectic groupoid $\Sigma$, whose germ is
unique up to isomorphism. The base $M$ of $\Sigma$ is a Poisson
submanifold and the Dirac structure on $N\subset M$ is the one induced
from $M$. 
\end{enumerate}

\begin{remark} The definition given here of a presymplectic groupoid, 
  does not guarantee that to each (integrable) Dirac structure there
  will be a unique presymplectic groupoid integrating it. To obtain
  uniqueness, one must require a non-degeneracy condition which can be
  expressed by saying that $\ker \omega$ intersects the fibers of $\s$
  and $\t$ transversely. After the submission of this work, the papers
  \cite{BCWZ,CX} have appeared. There one can find a serious
  discussion of presymplectic groupoids, multiplicative forms and their
  main properties.
\end{remark}

\subsection{Lie-Dirac submanifolds}
We consider now Poisson-Dirac submanifolds which admit a Dirac
projection compatible with the Lie brackets on
one-forms induced by the Poisson structures.

\begin{definition}
A submanifold $N\subset M$ of a Poisson manifold is called a
\textbf{Lie-Dirac submanifold} if it admits a Dirac projection $p:T_N
M\to TN$ such that $p^*:T^*N\to T^*M$ preserves Lie brackets (i.e., is
a Lie algebroid map). Such projections will be called
\textbf{Lie-Dirac projections}.
\end{definition}

Lie-Dirac submanifolds were first introduced by Xu in \cite{Xu1},
under the name ``Dirac submanifolds''. Note that while the property of
being a Poisson-Dirac submanifold is a local property, the property of being a
Lie-Dirac submanifold is a global property. The obstructions will be
studied below.

Let us first give alternative characterizations of these submanifolds. The
proof is immediate and is left to the reader.

\begin{proposition}
\label{prop:Lie:Dirac}
Let $M$ be a Poisson manifold and $N\subset M$ a submanifold. The
following statements are equivalent:
\begin{enumerate}
\item[(i)] $N$ is a Lie-Dirac submanifold;
\item[(ii)] There exists a Dirac projection $p:T_N M\to TN$ such that 
 $p([\Pi, X])= [\Pi_{N}, p(X)]$ for every $X\in\X(M)$;
\item[(iii)] There exists a bundle $E$ such that $T_N M=TN\oplus
  E$ and $E^0\subset T^*M$ is a Lie subalgebroid;
\item[(iv)] There exists a bundle $E$ such that $T_N M=TN\oplus
  E$ and $E$ is a coisotropic submanifold of $TM$.
\end{enumerate}
\end{proposition}

Note that for property (iv) one considers on $TM$ the tangent Poisson
structure.

\begin{example}
Any cosymplectic submanifold is a Lie-Dirac submanifold. However,
Poisson-Dirac submanifolds of constant rank, Poisson submanifolds, or even
symplectic leaves, may fail to be Lie-Dirac submanifolds.
\end{example}

A Poisson-Dirac submanifold $N\subset M$ of constant rank is
Lie-Dirac if and only if the sequence of Lie algebroids
(\ref{seq-Dirac}) has a Lie algebroid splitting. Since this splitting
integrates to a homomorphism of the corresponding Weinstein groupoids,
we can complete the diagram above to a commutative diagram
\[ 
\newdir{ (}{{}*!/-5pt/@^{(}}
\xymatrix{ \Sigma_N(M)\ar@{ (->}[r]\ar[dr]&\Sigma(M)\\& \Sigma(N)\ar[u]}\]
where the map going up is an embedding. 

Now note that even if $N$ does not have constant rank, the Lie-Dirac projection
induces a groupoid homomorphism $\Sigma(N)\to\Sigma(M)$ and we have
the following theorem, which is a slight improvement of Xu's results.

\begin{theorem} 
Let $M$ be an integrable Poisson manifold. Then any Lie-Dirac
submanifold $N\subset M$ is integrable, and $\Sigma(N)$ is a
symplectic subgroupoid of $\Sigma(M)$. More precisely, any Lie-Dirac
projection $p$ induces a groupoid embedding of $\Sigma(N)$ into
$\Sigma(M)$, which is a symplectomorphism onto a symplectic subgroupoid
of $\Sigma(M)$. Conversely, any such embedding is of this type.
\end{theorem}

\begin{proof}
The first part of the theorem is clear. So assume that 
$i:\Sigma'\hookrightarrow\Sigma(M)$ is a symplectic groupoid embedding
of a symplectic groupoid $\Sigma'$ over $N$.
We can identify $T^*N\simeq T_N\Sigma'/TN$ and $T^*M\simeq
T_M\Sigma/TM$. Hence, we obtain a Lie algebroid map
$i_*:T^*N\to T^*M$ whose image is a Lie subalgebroid $A\subset
T^*M$ which is transversal to $(TN)^0$. It is clear that $E=A^0$
satisfies (iii) of Theorem \ref{prop:Lie:Dirac}, so $N$ is a Lie-Dirac
submanifold and the embedding $i$ can be identified with the embedding
$\Sigma(N)\to \Sigma(M)$.
\end{proof}

Therefore, Lie-Dirac submanifolds of $M$ are the base manifolds of
symplectic subgroupoids of $\Sigma(M)$.  It should be noted that for a
given Lie-Dirac submanifold $N\subset M$ of a Poisson manifold there
can be \emph{distinct} connected symplectic subgroupoids of
$\Sigma(M)$ over $N$. In other words, the Poisson manifold $N$ does
not determine uniquely the subgroupoid of $\Sigma(M)$ (the image of
the embedding). This is, of course, because there can be several
distinct Lie-Dirac projections. Here is a very simple example.

\begin{example}
A manifold $M$ with the zero bracket integrates to the symplectic
groupoid $\Sigma=T^*M$, where $\omega$ is the canonical symplectic
structure, $\s=\t$ is the projection to $M$, and the group
operation is addition in the fibers. If $N\subset M$ is any
submanifold, then a vector subbundle $E\subset T^*M$ over $N$,
with $\rank E=\dim N$, is a symplectic subgroupoid
$\Sigma'\subset\Sigma$ if and only if $E\cap(TN)^\perp=\set{0}$. In this
case, $N$ is obviously a Poisson submanifold of $M$, but there are
many symplectic subgroupoids with the same base $N$.
\end{example}

As was already remarked by Xu in \cite{Xu1}, not every Poisson
submanifold is a Lie-Dirac submanifold. Here we give the obstruction
and show that it is related to the monodromy of Section
\ref{monodromy}. So let $i:N\hookrightarrow M$ be a Poisson
submanifold, giving rise to the short exact sequence of Lie algebroids
\begin{equation}
\label{eq:sequence:Poisson:submanifold} \xymatrix{0\ar[r]&
\nu^*(N)\ar[r]&T^*_N M\ar[r]^{i^*}& T^*N\ar[r]& 0.}
\end{equation}
Choosing a Dirac projection $p:T_N M\to TN$ is the same as choosing a
splitting $p^*:T^*N\to T^*_N M$ of this exact sequence.  This choice
gives rise to a contravariant connection $\nabla:\Omega^1(N)\times
\Gamma(\nu^*(N))\to \Gamma(\nu^*(N))$ defined by:
\[\nabla_{\al}\be =[p^*(\al),\be].\]
If the Lie algebra bundle $\nu^*(N)$ is abelian (this is the case,
for example, if $N$ consists only of regular points), this connection
is independent of the splitting, and it is flat:
\[ R_{\nabla}(\al,\beta)\equiv
\nabla_\al\nabla_\be-\nabla_\be\nabla_\al-\nabla_{[\al,\be]}=0.\]
In other words, $\nu^*(N)$ is canonically a flat Poisson vector
bundle over $N$. Then the curvature 2-form
\[ \Omega(\al,\beta)=p^*([\al,\be])-[p^*(\al),p^*(\be)],\]
defines a Poisson cohomology 2-class $[\Omega]\in
H^2_{\Pi}(N,\nu^*(N))$ which is the obstruction for $N$ to be a Lie-Dirac
submanifold. For example, for a regular leaf $L$ of a Poisson manifold
the Poisson cohomology $H^2_{\Pi}(L,\nu^*(L))$ coincides with the
cohomology $H^2(L,\nu^*(L))$ and we conclude that:

\begin{corollary}
Let $L$ be a regular leaf of a Poisson manifold. Then $L$ is a
Lie-Dirac submanifold if and only if the canonical cohomology class $[\Omega]\in
H^2(L,\nu^*(L))$ vanishes. If $L$ is simply connected, then $L$ is a
Lie-Dirac submanifold if and only if its monodromy group vanishes.
\end{corollary}

\begin{example}
Let us return to the Poisson brackets on $\mathfrak{su}^*(2)$ of
Section \ref{example:basic}. The origin is always a Lie-Dirac
submanifold. The spheres $S^2_R$, for $R>0$, are regular, compact
and simply connected leaves, and it follows that $S^2_R$ is a
Lie-Dirac submanifold if and only if the variation of the symplectic area $A'(R)$
vanishes.
\end{example}

Recall now (see e.g.~\cite{Vai}) that the
\textbf{characteristic form class} of the regular Poisson manifold
$M$ is the relative cohomology class $[\xi]=[d\omega]\in
H^3_{\text{rel}}(M,\F)$, where $\omega$ is the foliated symplectic form. The
characteristic form class is obviously the obstruction to the
existence of a closed 2-form on $M$ which pulls back to the
symplectic 2-form on each leaf. Equivalently, by the coisotropic
embedding theorem of Gotay \cite{Got}, it is the obstruction for the
existence of a leafwise symplectic embedding of $M$. From the remarks in
Section \ref{Integrability:regular}, we have a commutative diagram
relating the different classes associated with a regular Poisson
manifold
\[
\xymatrix{[\omega]\in H^2(\F) \ar[d]_{\delta}\ar[dr]^{d_{\nu}}\\
[\xi]\in H^3_{\text{rel}}(M,\F)\ar[r]&H^2(\F,\nu^*)\ni [\Omega]}
\]
where $\delta$ is defined from the long exact sequence of the pair $(M,\F)$:
\[
\xymatrix{\dots\ar[r]& H^k_{\text{rel}}(M,\F)\ar[r]& H^k(M)\ar[r]&
  H^k(\F)\ar[r]^\delta & H^{k+1}_{\text{rel}}(M,\F)\ar[r]&\dots}
\]
In particular we obtain:

\begin{corollary}
Let $M$ be a regular Poisson manifold. Then
\begin{enumerate}
\item[(i)] If $M$ admits a leafwise symplectic embedding then
every leaf of $M$ is a Lie-Dirac submanifold;
\item[(ii)] If every leaf of $M$ is a Lie-Dirac submanifold then $M$
is integrable.
\end{enumerate}
\end{corollary}

Note that the reverse implications in general do not hold.

%%%%%%%%%%%%%%%%%%%%%%%%%%%%%%%%%%%%%%%%%%%%%%%%%%%
%%%%%%%%%%%%%%%%%%%%%%%%%%%%%%%%%%%%%%%%%%%%%%%%%%%
\section{Morita equivalence}                      %
\label{Morita}                                   %
%%%%%%%%%%%%%%%%%%%%%%%%%%%%%%%%%%%%%%%%%%%%%%%%%%%
%%%%%%%%%%%%%%%%%%%%%%%%%%%%%%%%%%%%%%%%%%%%%%%%%%%
In this section we discuss Morita equivalence of Poisson
manifolds. This notion, introduced by P. Xu \cite{Xu2}, originates
on Weinstein's dual pairs as a global form of Lie's dual function
groups. Intuitively, Morita equivalence of Poisson manifolds means
isomorphism from the point of view of transversal Poisson geometry
(in particular, it induces Poisson diffeomorphisms between the
transversal Poisson structures). 

\subsection{Integrable Poisson manifolds}
Let us start by recalling Xu's definition:

\begin{definition}
\label{Morita:Xu} Two Poisson manifolds $M_1$ and $M_2$ are
\textbf{Morita equivalent} if there exists a two leg diagram
\[
\xymatrix{
&S\ar[dl]_{\pi_1}\ar[dr]^{\pi_2}\\
M_1& & \overline{M}_2}
\]
where $S$ is a symplectic manifold and the maps $\pi_1$ and
$\pi_2$ satisfy:
\begin{enumerate}
\item[(a)] each $\pi_i$ is a complete Poisson map; 
\item[(b)] each $\pi_i$ is a surjective submersion; 
\item[(c)] $\pi_1$ and $\pi_2$ have symplectic orthogonal, simply
  connected, fibers. 
\end{enumerate}
\end{definition}

From our section on symplectic realizations we see that, in a Morita equivalence, both 
Poisson manifolds have to be integrable, and, if $M$ is integrable, then the associated symplectic
groupoid $\Sigma(M)$ defines a Mosita equivalence of $M$ with itself ($\pi_{i}$ the source/the target). Hence

\begin{proposition} For a Poisson manifold $M$, the following are
equivalent:
\begin{enumerate}
\item[(i)] $M$ is Morita equivalent to another Poisson manifold;
\item[(ii)] $M$ is Morita equivalent to itself; 
\item[(iii)] $M$ is integrable.
\end{enumerate}
\end{proposition}

Moreover, if $S$ defines an equivalence between $M$ and $N$, we
also see that $S$ comes equipped with free (symplectic) actions of
$\Sigma(M_1)$ and $\Sigma(M_2)$, and the associated orbit spaces
are $M_1$, and $M_1$, respectively.
\[
\xymatrix{ \Sigma(M_1)\ar[d]\ar@<1ex>[d]
&S\ar[dl]_{\pi_1}\ar[dr]^{\pi_2}&
\Sigma(M_2)\ar[d]\ar@<1ex>[d]\\
M_1& & M_2}
\]
Or, in the terminology of \cite{Xu2,Xu3}, we recover Xu's result that
two integrable Poisson manifolds are Morita equivalent if and only
if their symplectic groupoids $\Sigma(M_1)$ and $\Sigma(M_2)$ are
symplectic Morita equivalent. Now we can proceed as in the ring-theoretic
version of Morita equivalence and compose such equivalences: if
$S_i$ defines an equivalence between $M_{i}$ and $M_{i+1}$, $i\in
\{ 1, 2\}$, then $S_{1}\otimes_{\Sigma(M_2)} S_{2}$, that is, the
quotient of $S_{1}\times_{M_{2}}S_{2}$ by the diagonal action of
$\Sigma(M_{2})$, is an equivalence between $M_1$ and $M_3$. In
particular:

\begin{corollary} On the class of integrable Poisson manifolds,
Morita equivalence is an equivalence relation.
\end{corollary}

\subsection{Weak Morita Equivalence}
{For} general Poisson manifolds, there are two important questions one
should address: 
\begin{itemize} 
\item Find a satisfactory notion of Morita equivalence which works
well also for non-integrable Poisson manifolds. 
\item Construct ``Morita invariants'', i.e., invariants which allow us to
distinguish non-Morita equivalent Poisson manifolds and eventually
classify Morita equivalence classes.
\end{itemize}

The notion of weak Morita equivalence of Poisson manifolds is very
similar to that of Morita equivalence of foliations (discussed in
\cite{Cra, Ginz}) which, in turn, is an infinitesimal version of Morita
equivalence of holonomy groupoids (see \cite{Haefl}).  More precisely,
given a submersion $\phi: Q\to M$ and a (regular) foliation $\F$ on
$M$, the pull-back foliation on $Q$, denoted $\phi^{\star}\F$, is the
foliation whose leaves are the connected components of $\phi^{-1}(L)$,
with $L$ leaf of $\F$. In terms of the associated involutive
subbundles, $\phi^{\star}\F$ consists of vectors $X$ tangent to $Q$
with the property that $d\phi(X)\in \F$, and is in general different
from the pull-back vector bundle $\phi^*\F$ (note the difference in
the notation).  We will say that two (regular) foliations $\F_{i}$ on
$M_{i}$ ($i\in \{1, 2\}$) are Morita equivalent if there exists a
manifold $Q$, and submersions $\pi_{i}$ from $Q$ onto $M_{i}$ with
simply connected fibers: 
\[ 
M_1\stackrel{\pi_1}{\longleftarrow}Q\stackrel{\pi_2}\longrightarrow M_2 , 
\] 
and such that $\pi_{1}^{*}\F=\pi_{2}^{*}\F$.

For Poisson manifolds we proceed similarly. First of all, for a
submersion $\phi: Q\to M$ from a manifold $Q$ into a Poisson manifold
$M$, we form the bundle $\phi^\star T^*M$ over $Q$:
\[ 
\phi^\star T^*M=\set{(\al,X)\in \phi^*T^*M\times TQ: \#\al=\phi_*X}.
\]
This is a Lie algebroid over $Q$ with the anchor $(\al,X)\mapsto
X$, and with the Lie bracket
\[
[(f\phi^*\al,X),(g\phi^*\be,Y)]=
(fg\phi^*[\al,\be]+X(g)\phi^*\be-Y(f)\phi^*\al,[X,Y]),
\]
where $f,g\in C^\infty(Q)$, $X,Y\in\X(Q)$, $\al,\be\in \Omega^1(M)$,
and $\phi^*\al$ is the induced section of $\phi^*T^*M$. The idea is
that this should represent the pull-back of the Poisson structure to
$Q$. The point is that pull backs of Poisson structures generally only
make sense as Dirac structures, and $\phi^\star T^*M$ is precisely
this pull-back Dirac structure, viewed as a Lie algebroid.

% First of all, for a submersion $\phi: Q\to M$ from a manifold
% $Q$ into a Poisson manifold $M$, one forms the \emph{pull-back Lie
% algebroid} $\phi^\star T^*M$ as follows (note the difference with
% the pull-back vector bundle):
% \[ \phi^\star T^*M=\set{(\al,X)\in \phi^*T^*M\times TQ:
%   \#\al=\phi_*X},\]
% with the anchor $(\al,X)\mapsto X$, and with the Lie bracket
%   \[ [(f\al,X),(g\be,Y)]=(fg[\al,\be]+X(g)\be-Y(f)\al,[X,Y]),\]
% where $f,g\in C^\infty(Q)$, $X,Y\in\X(Q)$ and $\al,\be$ are
% pull-backs of 1-forms on $M$. Alternatively, $\phi^\star T^*M$ is
% the algebroid underlying the pull-back Dirac structure on $Q$.

\begin{definition}
Let $M_1$ and $M_2$ be two Poisson manifolds. We say that they are
weakly 
Morita equivalent if there exists a manifold $Q$ together
with two submersions onto $M_1$ and $M_2$ with simply connected
fibers
\[ 
M_1\stackrel{\pi_1}{\longleftarrow}Q\stackrel{\pi_2}\longrightarrow M_2 , 
\] 
such that the pull-back algebroids (equivalentely, Dirac structures)
on $P$ are isomorphic: 
\[ \pi_1^\star T^*M_1\simeq \pi_2^\star T^*M_2.\]
\end{definition}
The reader will notice that Dirac structures furnish a natural setting
fot Morita equivalence. We refer to \cite{BuRa} for details on this approach.
Also, the previous notion of Morita equivalence originates in \cite{Cra}
(see Theorem 2 there and the comment preceeding it) and it coincides with 
the equivalence used in \cite{Ginz}.

\begin{proposition}
\label{iso:pull-backs}
 \text{}
\begin{enumerate}
\item[(i)] Morita equivalence implies weak Morita equivalence;
\item[(ii)] weak Morita equivalence is an equivalence relation on
the class of all Poisson manifolds.
\end{enumerate}
\end{proposition}

\begin{proof}
Let $S$ be as in the definition of Morita equivalence. We denote
by $\rho_{i}: \pi^*T^*M_{i}\to TS$ the induced actions on $S$.
There are bundle isomorphisms
\[ \pi_1^\star T^*M_1
\simeq \pi_1^*T^*M_1 \oplus \pi_2^*T^*M_2 \simeq \pi_2^\star
T^*M_2, \] which, for sections $\al$ of $\pi_1^*T^*M_1$ and $\be$
of $\pi_2^*T^*M_2$ which are pull-backs of one forms on $M_1$ and
$M_2$, are given by
\[ (\al,\rho_1(\al)+\rho_2(\be))\longmapsfrom (\al,\be)\longmapsto
(\be,\rho_1(\al)+\rho_2(\be)).\] A careful computation shows that
this is actually a Lie algebroid isomorphism. Assume now that $P$
defines another weak Morita equivalence between $M_{2}$ and
$M_{3}$. We then form the fibered product $R= Q\times_{M_2}P$ and
denote by $q$ and $p$, the projections into $Q$, and $P$,
respectively:
\[ \xymatrix{
 & & R\ar[ld]_{q}\ar[rd]^{p} & & \\
 & Q\ar[ld]_{\pi_1}\ar[rd]^{\pi_2} & & P\ar[ld]_{\pi_3}\ar[rd]^{\pi_4} & \\
M_1 & & M_2 & & M_{3} }
\]
Then $R$, together with $\pi_1 q$ and $\pi_4 p$ defines a weak
Morita equivalence. This follows from the functoriality of the
operation $\phi^\star$ on Lie algebroids.
\end{proof}

We now proceed with the description of various (weak) Morita
invariants. Due to the similarity with the notion of Morita equivalence
of foliations, it is not surprising that:

\begin{proposition}
Let $M_1$ and $M_2$ be weakly Morita equivalent Poisson manifolds with
symplectic foliations $\F_1$ and $\F_2$. Then any weak Morita
equivalence $Q$ between them induces a homeomorphism
\[ \phi_{Q}:M_1/\F_1\simeq M_2/\F_2,\]
and the fundamental groups of corresponding leaves are isomorphic
\[  \phi_Q(L_1)=L_2 \Longrightarrow \pi(L_{1})\cong \pi(L_{2}).\]
Moreover, if one of the two leaf spaces
is a smooth manifold, then
so is the other one, and the map above is a diffeomorphism. 
\end{proposition}

\begin{proof}
It is easy to see that given a submersion $\phi:Q\to M$ onto a
Poisson manifold $M$, the leaves of $\phi^\star T^*M$ are of the
form $\phi^{-1}(L)$, where $L$ is a leaf of $M$. Hence, given a
Morita equivalence
$M_1\stackrel{\pi_1}{\longleftarrow}Q\stackrel{\pi_2}{\longrightarrow}M_2$,
the map $L\mapsto \pi_2(\pi_1^{-1}(L))$ is a bijective
correspondence between the leaves of $M_1$ and the leaves of
$M_2$. This bijection makes the following diagram commute:
\[
\xymatrix{
&Q\ar[dl]_{\pi_1}\ar[dr]^{\pi_2}\\
M_1\ar[d]& & M_2\ar[d]\\
M_1/\F_1\ar@{<->}[rr]& &M_2/\F_2 }
\]
The statements on the leaf spaces follow from the fact that the maps
going down are quotient maps. On the other hand, the statement about
the fundamental groups of corresponding leaves, follows from the fact
that any submersion with simply connected fibers induces isomorphisms
in the first homotopy groups (cf. the Appendix in \cite{GiGo}).
\end{proof}

Our next result shows that a large number of geometric invariants are
weak Morita invariants, including the monodromy groups we have
introduced before. In this respect, note that the higher homotopy groups of
corresponding leaves may not be isomorphic, so one may think of the
monodromy groups (which come from the second homotopy groups of the
leaves), as the next level in the hierarchy. 

\begin{theorem}
\label{monodromy-Morita} 
The monodromy groups, isotropy Lie algebras and groups, first Poisson
cohomology groups and integration along cotangent 
paths, are all weak
Morita invariants. More precisely, let $Q$ be a weak Morita
equivalence between $M_1$ and $M_2$, and let
$x\in M_1$, $y\in M_2$ be points whose
leaves correspond. Then the following are isomorphic:
\begin{enumerate}
\item[(i)] The monodromy groups $\NN_{x}(M_1)$ and $\NN_{y}(M_2)$;
\item[(ii)] the isotropy Lie algebras $\mathfrak{g}_{x}$, and
$\mathfrak{g}_{y}$; 
\item[(iii)] the isotropy groups $\Sigma(M_1, x)$ and
$\Sigma(M_2, x)$; 
\item[(iv)] the reduced isotropy groups $\Sigma^0(M_1,x)$ and
  $\Sigma^0(M_2, x)$;
\item[(v)] the Poisson cohomology groups $H^{1}_{\Pi}(M_1)$ and
  $H^{1}_{\Pi}(M_2)$, and the two maps 
\[ \Sigma(M_i, x_i)\times H^{1}_{\Pi}(M_i) \rmap \mathbb{R}, ([a], X)\mapsto
\int_{a}X \quad (i=1,2).\]
\end{enumerate}
\end{theorem}

\begin{proof}
The idea is the same as in the previous proof:
given a Lie algebroid $A$ over $M$, and a submersion $\phi: Q\to M$
with simply connected fibers, we prove that all these groups for
$\phi^\star A$ are isomorphic to the ones of $A$. We use here that all
we have said about $T^*M$ and $\Sigma(M)$ works for general Lie
algebroids (cf.~\cite{CrFe}). 

We have already remarked that the leaves of $\phi^\star A$ are of type
$\phi^{-1}L$, with $L$ leaf of $A$. Since the fibers are simply
connected, the maps $\pi_{2}(\phi^{-1}L)\to \pi_{2}(L)$ are
surjective \cite{Ginz}. 

Next, the kernel of the anchor of $\phi^\star A$ is clearly isomorphic
to the kernel of the anchor of $A$, i.e., the isotropy Lie algebras
$\mathfrak{g}_{q}(\phi^\star A)$ and $\mathfrak{g}_{\phi(q)}(A)$ are
isomorphic. For the other groups, due to the canonical sequences
(\ref{eq1}), (\ref{eq2}) (and their obvious versions for Lie
algebroids), and the fact that $\phi$ induces isomorphism in the first
homotopy groups, it is enough to prove the statement for the reduced
isotropy group $\Sigma^0(M, x)$ (and its Lie algebroid versions,
denoted $G^{0}(A, x)$). Both groups $G(\phi^\star A, q)$ and $G(A, x)$
($x= \phi(q)$) are quotients of the same space (of paths in the
isotropy Lie algebra), and one only has to check that the equivalence
relation (homotopy) is the same. For this one uses that, when working
in the smooth category, submersions do behave like Serre fibrations
(cf.~the Appendix in \cite{GiGo}). Alternatively, one can prove all
these isomorphisms at once, computing the groupoid $G(\phi^{\star}A)$
of $\phi^{\star}A$ (see \cite{CrFe} for notations): 

First of all, a path for $\phi^{\star}A$ is a pair $(a,
\tilde{\gamma})$, where $A$ is an $A$-path, and $\tilde{\gamma}$ is a
path in $Q$ over the base path of $a$.  However, since
\begin{itemize}
\item[(a)] for any two paths $\tilde{\gamma}_{0}$ and $\tilde{\gamma}_{1}$
in $Q$ covering the same path $\gamma$ in $M$, there is a homotopy
$\tilde{\gamma}_{\epsilon}$ between them, covering $\gamma$ (this
follows since the fibers are simply connected, see the Appendix in
\cite{Ginz});
\item[(b)] if $\tilde{\gamma}_{\epsilon}$ is as above, and $a$ is an
$A$-path with base path $\gamma$, then $(a,
\tilde{\gamma}_{\epsilon})$ is a $\phi^{\star}A$-homotopy.
\end{itemize}
It follows that, when looking at homotopy classes of
$\phi^{\star}A$-paths $(a, \tilde{\gamma})$, it is only $a$ and
the end points of $\tilde{\gamma}$ that matter. Then, working with
triples $(p, a, q)$, where $a$ is an $A$-path, and $\phi(p)=
\gamma(0)$, $\phi(q)= \gamma(1)$, and writing down what
$\phi^{\star}A$-homotopy means, ones gets $G(\phi^{\star}A)=
Q\times_{M}G(A)\times_{M}Q$ consists of triples $(p, g, q)$ with
$\phi(p)= \s (g)$, $\phi(q)= \t (g)$, with the multiplication $(q,
h, r) (p, g, q)= (p, hg, r)$. This immediately imply the rest.

Finally, the statement about Poisson cohomology follows from the fact
(see \cite{Cra}, Theorem 2) that if $\phi:Q\to M$ is a submersive
submersion with simply-connected fibers then $A$ and $\phi^\star A$
have isomorphic first cohomology groups.
\end{proof}

\begin{example} Let us look at some examples:
\begin{enumerate}
\item[1.] The fundamental group $\pi_1(S)$ is a complete Morita
invariant in the symplectic case. 
\item[2.] The isotropy Lie algebra at zero is a complete invariant in the
case of linear structures (we actually see that $\mathfrak{g}^*$ and
$\mathfrak{h}^*$ are weakly Morita equivalent around the
origin if and only if and
only if $\mathfrak{g}$ and $\mathfrak{h}$ are isomorphic Lie
algebras).
\item[3.] Consider the sphere $S^2$ with its standard symplectic structure
(area form) and form the Poisson-Heisenberg manifold $M(S^2)$ (see the
example in \ref{Heisenberg}). This Poisson manifold is not Morita
equivalent to the any of the Poisson manifolds $M_a$ considered in the
example of Section \ref{example:basic}, just because they have
different leaf spaces.
\item[4.] The Poisson manifolds $M_{a}$, for different $a$'s, all have the
same leaf spaces. However, one can find $a$'s for which the monodromy
groups vary differently, hence the associated Poisson structures are
not weakly Morita equivalent (this will actually force one of the
manifolds to be non-integrable, but this can be avoided by
going to leaf spaces of
higher dimensions). 
\end{enumerate}
\end{example}

We say that a Poisson manifold $M$ has simple symplectic foliation if
$\F$ is regular, and the leaf space $B= M/\F$ is smooth. Equivalently,
one has a submersion $\pi: M\to B$ with connected symplectic fibers,
and $M$ has the induced Poisson structure. Note that in this case
there are canonical identifications $\nu_{x}\cong T_{\pi(x)}B$, and
the monodromy groups $\NN_{x}$ define a bundle of subgroups
\[ \NN(M)\subset T^{*}B .\]
We know that the continuity of $\NN(M)$ is related to the
integrability of $M$. It is not difficult to see that the class of
such Poisson manifolds is closed under weak Morita equivalence.
{From} the previous proofs we see that

\begin{corollary}
\label{simple-foliations} For Poisson manifolds $M$ with simple
symplectic foliation, the pair $(B, \NN(M))$ is a weak Morita
invariant. More precisely, a weak Morita equivalence $Q$ between
two such $M_{1}$ and $M_{2}$ induces a diffeomorphism $\phi:
B_{1}\to B_{2}$ with $(d\phi)^{*}_{x}(\NN(M_{1}))= \NN(M_{1})$.
\end{corollary}

Specializing to product foliations, Xu proved a similar result for ordinary Morita
  equivalence, along with a converse (see also the next proposition). 

\begin{remark} It is important that the isomorphism between the monodromy
is induced by the differential. For instance, consider the Poisson 
manifolds $M_{a}^{'}= M_{a}\setminus \{ 0\}$. They all have the same
leaf space, and one can choose $a$ and $b$ such that the corresponding 
Poisson manifolds have the same (i.e. isomorphic)  monodromy, without
being Morita equivalent.
\end{remark}

The following proposition suggests that, in some cases, the leaf space, the
homotopy groups of the leaves, and the monodromy groups form a
\emph{complete} set of weak Morita invariants

\begin{proposition}
\label{complete-invariant} 
Let $M$ be a Poisson manifold with simple symplectic foliation $\F$
and compact simply connected leaves. If $\NN(M)= 0$, then $M$ with the
induced Poisson structure is weakly Morita equivalent to $B=M/\F$ with
the zero Poisson structure.
\end{proposition}

\begin{proof} 
Let $\sigma: \F\to T^*M$ be a splitting of the anchor, and let
$\Omega\in \Omega^2(\F; \nu^*)$ be its curvature.  Since each leaf $L$
is simply connected, the condition on the monodromy shows that
$\Omega|_{L}$ is exact. From the Reeb stability theorem each $L$ has a
neigborhood of type $T\times L$, hence $\Omega$ is exact in a
neigborhood $L$. By choosing a partition of unity supported on such
opens neigborhoods, and which is constant on each leaf (e.g.  the
pullback by $\pi:M\to B$ of an open cover of the base), we see that
$\Omega$ is exact. Hence $\sigma$ has a global splitting compatible
with the brackets. This shows that $\Sigma(M)= \pi(\F)\times_{M}
\nu^*$, which is clearly Morita equivalent to $T^*B$ (as a bundle of
abelian Lie groups over $B$).
\end{proof}

\begin{remark} Our aim here was not only to describe Morita
invariants, but also to point out that all the algebraic
invariants we know are weak Morita invariants. Nevertheless,
although this notion does behave well for all Poisson manifolds,
it is not our intention to present it as the ``satisfactory''
notion of Morita equivalence that we asked at the beginning of this
section, for it does not take the symplectic/Poisson picture fully into account.
Let us point out two other possible notions of Morita equivalence for
non-integrable Poisson manifolds:
\begin{enumerate}
\item[] \textbf{Symplectic Morita Equivalence:} We mimic Xu's
  definition (see Definition \ref{Morita:Xu}) but we remove the
  completeness assumption. Then we obtain a symmetric relation which
  is also reflexive (because of the existence of a local symplectic
  groupoid integrating any $M$). It is however not clear that this
  relation is transitive,
so we should take the equivalence relation it
  generates.
\item[] \textbf{Algebraic Morita Equivalence:} For $M_1$ and $M_2$ to
  be algebraic Morita equivalent we require the existence of a Lie
  algebroid bi-module for $T^*M_1$ and $T^*M_2$: there are commuting
  free actions $\rho_1$ and $\rho_2$ of $T^* M_1$ and $T^* M_2$ on
  $Q$, with moment maps
  $M_1\stackrel{\pi_1}{\longleftarrow}Q\stackrel{\pi_2}{\longrightarrow}M_2$,
  where each $\pi_i$ is a surjective submersion with simply connected
  fibers, such that the orbits of $\rho_1$ (respectively, $\rho_2$) are
  the fibers of $\pi_2$ (respectively, $\pi_1$).
\end{enumerate}

It is not hard to see that one has the following chain of
implications:
\[ 
\text{symplectic Morita }\Longrightarrow
\text{ algebraic Morita }\Longrightarrow
\text{ weak Morita }
\]
In our opinion, it is an important open question to show
that these relations coincide or else to find new Morita invariants
which are not weak Morita invariants.  
\end{remark}
% -----------------------------------------------------------------------
\bibliographystyle{amsplain}
\def\lllll{}

\end{document}